\input lanlmac
\def\href#1#2{{#2}}

\input epsf.tex

\overfullrule=0mm

\newcount\figno
\figno=0
\def\fig#1#2#3{
\par\begingroup\parindent=0pt\leftskip=1cm\rightskip=1cm\parindent=0pt
\baselineskip=11pt
\global\advance\figno by 1
\midinsert
\epsfxsize=#3
\centerline{\epsfbox{#2}}
\vskip 12pt
{\bf Fig.\ \the\figno:} #1\par
\endinsert\endgroup\par
}
\def\figlabel#1{\xdef#1{\the\figno}}
\def\encadremath#1{\vbox{\hrule\hbox{\vrule\kern8pt\vbox{\kern8pt
\hbox{$\displaystyle #1$}\kern8pt}
\kern8pt\vrule}\hrule}}


\def\IR{\relax{\rm I\kern-.18em R}}
\font\cmss=cmss10 \font\cmsss=cmss10 at 7pt

\font\cmss=cmss10 \font\cmsss=cmss10 at 7pt
\def\IZ{\relax\ifmmode\mathchoice
{\hbox{\cmss Z\kern-.4em Z}}{\hbox{\cmss Z\kern-.4em Z}}
{\lower.9pt\hbox{\cmsss Z\kern-.4em Z}}
{\lower1.2pt\hbox{\cmsss Z\kern-.4em Z}}\else{\cmss Z\kern-.4em Z}\fi}
\def\IN{\relax{\rm I\kern-.18em N}}
\def\b{\circ}
\def\n{\bullet}

\def\gbbbb{\Gamma_4^{\hbox{$\scriptstyle \b \b$}\kern -8.2pt
\raise -4pt \hbox{$\scriptstyle \b \b$}}}
\def\gnnnn{\Gamma_4^{\hbox{$\scriptstyle \n \n$}\kern -8.2pt  
\raise -4pt \hbox{$\scriptstyle \n \n$}}}
\def\gnnnnnn{\Gamma_6^{\hbox{$\scriptstyle \n \n \n$}\kern -12.3pt
\raise -4pt \hbox{$\scriptstyle \n \n \n$}}}
\def\gbbbbbb{\Gamma_6^{\hbox{$\scriptstyle \b \b \b$}\kern -12.3pt
\raise -4pt \hbox{$\scriptstyle \b \b \b$}}}
\def\gbbbbc{\Gamma_{4\, c}^{\hbox{$\scriptstyle \b \b$}\kern -8.2pt
\raise -4pt \hbox{$\scriptstyle \b \b$}}}
\def\gnnnnc{\Gamma_{4\, c}^{\hbox{$\scriptstyle \n \n$}\kern -8.2pt
\raise -4pt \hbox{$\scriptstyle \n \n$}}}
\def\Rbud#1{{\cal R}_{#1}^{-\kern-1.5pt\blacktriangleright}}
\def\Rleaf#1{{\cal R}_{#1}^{-\kern-1.5pt\vartriangleright}}
\def\Rbudb#1{{\cal R}_{#1}^{\circ\kern-1.5pt-\kern-1.5pt\blacktriangleright}}
\def\Rleafb#1{{\cal R}_{#1}^{\circ\kern-1.5pt-\kern-1.5pt\vartriangleright}}
\def\Rbudn#1{{\cal R}_{#1}^{\bullet\kern-1.5pt-\kern-1.5pt\blacktriangleright}}
\def\Rleafn#1{{\cal R}_{#1}^{\bullet\kern-1.5pt-\kern-1.5pt\vartriangleright}}
\def\Wleaf#1{{\cal W}_{#1}^{-\kern-1.5pt\vartriangleright}}
\def\Cleaf{{\cal C}^{-\kern-1.5pt\vartriangleright}}
\def\Cbud{{\cal C}^{-\kern-1.5pt\blacktriangleright}}
\def\Crleaf{{\cal C}^{-\kern-1.5pt\circledR}}


\magnification=\magstep1
\baselineskip=12pt
\hsize=6.3truein
\vsize=8.7truein
 at 8truept
 at 8truept
 at 10truept

\font\bigrm=cmr12 at 14pt
\centerline{\bigrm Blocked edges on Eulerian maps and mobiles:}
\centerline{\bigrm Application to spanning trees,}
\centerline{\bigrm hard particles and the Ising model}

\bigskip\bigskip

\centerline{J. Bouttier, P. Di Francesco and E. Guitter}
  \smallskip
  \centerline{Service de Physique Th\'eorique, CEA/DSM/SPhT}
  \centerline{Unit\'e de recherche associ\'ee au CNRS}
  \centerline{CEA/Saclay}
  \centerline{91191 Gif sur Yvette Cedex, France}
\centerline{\tt bouttier@spht.saclay.cea.fr}
\centerline{\tt philippe@spht.saclay.cea.fr}
\centerline{\tt guitter@spht.saclay.cea.fr}

  \bigskip


     \bigskip\bigskip

     \centerline{\bf Abstract}
     \smallskip
     {\narrower\noindent
We introduce Eulerian maps with blocked edges as a general way to implement
statistical matter models on random maps by a modification of intrinsic
distances. We show how to code these dressed maps by means of mobiles, i.e. decorated trees 
with labeled vertices, leading to a closed system of recursion relations for 
their generating functions. We discuss particular solvable cases in detail, 
as well as various applications of our method to several statistical systems
such as spanning trees on quadrangulations, mutually excluding particles on 
Eulerian triangulations or the Ising model on quadrangulations.
\par}

     \bigskip

\nref\TUT{W. Tutte,
{\it A Census of planar triangulations} Canad. J. of Math. {\bf 14} (1962) 21-38;
{\it A Census of Hamiltonian polygons} Canad. J. of Math. {\bf 14} (1962) 402-417;
{\it A Census of slicings}, Canad. J. of Math. {\bf 14} (1962) 708-722;
{\it A Census of Planar Maps}, Canad. J. of Math. {\bf 15} (1963) 249-271.
}
\nref\QGRA{V. Kazakov, {\it Bilocal regularization of models of random
surfaces}, Phys. Lett. {\bf B150} (1985) 282-284; F. David, {\it Planar
diagrams, two-dimensional lattice gravity and surface models},
Nucl. Phys. {\bf B257} (1985) 45-58; J. Ambjorn, B. Durhuus and J. Fr\"ohlich,
{\it Diseases of triangulated random surface models and possible cures},
Nucl. Phys. {\bf B257}(1985) 433-449; V. Kazakov, I. Kostov and A. Migdal
{\it Critical properties of randomly triangulated planar random surfaces},
Phys. Lett. {\bf B157} (1985) 295-300.}
\nref\DGZ{P. Di Francesco, P. Ginsparg and J. Zinn--Justin, {\it 2D Gravity and Random Matrices},
Physics Reports {\bf 254} (1995) 1-131.}
\nref\BIPZ{E. Br\'ezin, C. Itzykson, G. Parisi and J.-B. Zuber, {\it Planar
Diagrams}, Comm. Math. Phys. {\bf 59} (1978) 35-51.}
\nref\CORV{R. Cori and B. Vauquelin, {\it Planar maps are well labeled trees},
Canad. J. Math. {\bf 33 (5)} (1981) 1023-1042.}
\nref\ARQUES{D. Arqu\`es, {\it Les hypercartes planaires sont des arbres 
tr\`es bien \'etiquet\'es}, Discr. Math. {\bf 58}(1) (1986) 11-24.}  
\nref\SCH{G. Schaeffer, {\it Bijective census and random
generation of Eulerian planar maps}, Electronic
Journal of Combinatorics, vol. {\bf 4} (1997) R20; see also
{\it Conjugaison d'arbres
et cartes combinatoires al\'eatoires}, PhD Thesis, Universit\'e 
Bordeaux I (1998).}
\nref\CONST{M. Bousquet-M\'elou and G. Schaeffer,
{\it Enumeration of planar constellations}, Adv. in Applied Math.,
{\bf 24} (2000) 337-368.}
\nref\CENSUS{J. Bouttier, P. Di Francesco and E. Guitter, {\it Census of planar
maps: from the one-matrix model solution to a combinatorial proof},
Nucl. Phys. {\bf B645}[PM] (2002) 477-499, arXiv:cond-mat/0207682.}
\nref\CS{P. Chassaing and G. Schaeffer, {\it Random Planar Lattices and 
Integrated SuperBrownian Excursion}, 
Probability Theory and Related Fields {\bf 128(2)} (2004) 161-212, 
arXiv:math.CO/0205226.}
\nref\MOB{J. Bouttier, P. Di Francesco and E. Guitter. {\it Planar maps as labeled mobiles},
Elec. Jour. of Combinatorics {\bf 11} (2004) R69, arXiv:math.CO/0405099.}
\nref\BMS{M. Bousquet-M\'elou and G. Schaeffer,{\it The degree distribution
in bipartite planar maps: application to the Ising model},
arXiv:math.CO/0211070.}
\nref\HPCOM{J. Bouttier, P. Di Francesco and E. Guitter. {\it Combinatorics of Hard Particles
on Planar Graphs}, Nucl.Phys. {\bf B655} (2003) 313-341, arXiv:cond-mat/0211168.}
\nref\HObipar{J. Bouttier, P. Di Francesco and E. Guitter. {\it Combinatorics of bicubic maps 
with hard particles}, J.Phys. A: Math.Gen. {\bf 38} (2005) 4529-4560, arXiv:math.CO/0501344.}
\nref\HPMAT{J. Bouttier, P. Di Francesco and E. Guitter. {\it Critical and Tricritical Hard 
Objects on Bicolorable Random Lattices: Exact Solutions}, J.Phys. A: Math.Gen. {\bf 35} 
(2002) 3821-3854, arXiv: cond-mat/0201213.}
\nref\DK{B. Duplantier and I.K. Kostov, {\it Conformal spectra of polymers on a random surface},
Phys. Rev. Lett. {\bf 61} (1988) 1433-1437;  {\it Geometrical critical phenomena on a random 
surface on arbitrary genus}, Nucl. Phys. {\bf B340} (1990) 491-541.}
\nref\BOUKA{D. Boulatov and V. Kazakov, {\it The Ising model
on a random planar lattice: the structure of the phase 
transition and the exact critical exponents}, Phys. Lett. {\bf B186} (1987)
379-384.}
\nref\ON{
I.K. Kostov,  {\it O(n) vector model on a planar random lattice: spectrum of 
anomalous dimensions}, Mod.Phys.Lett. {\bf A4} (1989) 217-226; 
B. Eynard and J. Zinn-Justin, {\it The O(n) model on a random surface: critical 
points and large order behaviour}, Nucl. Phys. {\bf B386} (1992) 558-591, arXiv:hep-th/9204082;
B. Eynard, C. Kristjansen, {\it Exact solution of the O(n) model on 
a random lattice} , Nucl. Phys. {\bf B455} (1995) 577-618, arXiv:hep-th/9506193.}
\nref\POTTS{J.M. Daul, {\it  Q-states Potts model on a random planar lattice}, 
arXiv:hep-th/9502014;
G. Bonnet and B. Eynard, {\it The Potts-q random matrix model : loop equations, 
critical exponents, and rational case}, Phys.Lett. {\bf B463} (1999) 273-279, 
arXiv:hep-th/9906130.}
\nref\SIXV{
V. Kazakov and P. Zinn-Justin, {\it Two-Matrix model with ABAB interaction}, Nucl. Phys. {\bf B546}
(1999) 647-668, arXiv:hep-th/9808043; 
I.K. Kostov, {\it Exact Solution of the Six-Vertex Model on a Random Lattice}, 
Nucl. Phys. {\bf B575} (2000) 513-534, arXiv:hep-th/9911023.}
\newsec{Introduction}

Enumeration of maps, i.e. connected graphs embedded in a surface, is a fundamental issue 
that has drawn the attention of mathematicians and physicists all over the years ever since 
the seminal series of papers of Tutte on the subject \TUT. In physical applications, the subject
was extended so as to include {\it statistical models} defined on the maps, 
by equipping them with ``matter" configurations, say of particles, dimers, spins or trees.
These dressed maps were introduced as discrete realizations of random surfaces
in the general framework of two-dimensional quantum gravity (2DQG) \QGRA\ for which many
results have been obtained so far (see \DGZ\ for a review).
Several techniques have been used to address this problem of counting maps, 
with their possible matter decorations generally weighted by Boltzmann factors.   
Beside the original recursive decomposition technique of Tutte, mainly applied
to maps without matter, the most powerful approach was undoubtedly that based
on matrix integrals \BIPZ. Many results obtained that way display a simplicity
contrasting with the high sophistication of the method. In particular, most
planar enumeration results involve algebraic generating functions. This led combinatorists
to look for alternative and more constructive approaches that could explain this
simplicity [\xref\CORV,\xref\ARQUES]. Quite recently, a general technique was developed that relies
on bijections between planar maps with possible decorations and families of 
decorated trees \SCH. The bijections are of two different types. A first
class of bijections [\xref\SCH-\xref\CENSUS], well adapted to maps with prescribed vertex 
valences, consists
of a cutting algorithm of the map into a {\it blossom tree} that carries the minimal
information needed to recover the map. A second class of bijections [\xref\CS,\xref\MOB] 
is well adapted 
to the dual version of the former, namely maps with prescribed face valences. 
It makes use of the intrinsic geodesic distance on the map to code it by a so-called 
{\it labeled mobile}, i.e. a tree decorated by integer labels recording this distance 
in a way that allows to reconstruct the original map. 
Both methods were used to enumerate families of bicolored planar maps, or equivalently
their dual Eulerian maps, {\it in the absence of matter}, all displaying the
standard ``square-root singularity" characteristic of pure gravity (if we omit non-generic,
non-physical multicritical points). Still, one may in some cases 
reinterpret the color as a matter degree of freedom via suitable weightings. 
This gives access for instance to the solution of the Ising model on tetravalent maps \BMS, 
leading to a new singular behavior of the generating function at the price of 
some analytic continuation to reach the transition point. Similarly, this allows to treat 
the case of hard particles on arbitrary (non bicolored) maps, but with no new singular 
behavior in this case \HPCOM. 

In Ref. \HObipar, it was found how to extend the bijection with blossom trees to the case 
of bicolored maps equipped with interacting particles, allowing in particular to
recover in a purely combinatorial way the crystallization transition of hard particles.
Considering more involved particle exclusion rules leads to most of the possible 
universality classes of rational conformal theories (RCFT) with central charge $c<1$ coupled
to 2DQG \HPMAT. 

The aim of this paper is to extend the mobile formalism to this same family of
bicolored (or dually Eulerian) maps with mutually excluding particles. 
As we shall see however, our method
proves more general and applies to other classes of matter such as spanning
trees or Ising spins, here without analytic continuation, thus providing
a unified framework for all these models.
More precisely, as already mentioned, the original mobile coding of Eulerian
maps uses the geodesic distance from a given origin vertex on the map. 
A way to incorporate matter is precisely to modify this distance by introducing
{\it blocked edges}, i.e. edges excluded from geodesic paths. All the matter
systems above can indeed be reformulated in terms of configurations of blocked
edges on Eulerian maps. The mobile construction is then easily modified so
as to account for the presence of the blocked edges and this eventually leads to an
efficient coding of the maps at hand in the form of labeled trees with marked edges.
Using the tree structure, one can derive recursive equations for the corresponding generating
functions, leading to a combinatorial solution of the original problems.

The paper is organized as follows. Section 2 is devoted to the general theory
of map coding via mobiles. We first define precisely in section 2.1 the notion of
Eulerian maps with blocked edges and show in section 2.2 how to code them by 
labeled mobiles that we characterize in detail. The procedure leading back from
the mobile to the original map is described in section 2.3. In section 2.4, 
we further use the characterization of mobiles to obtain a system of recursive equations 
for their generating functions, also easily translated into map generating functions.

Sections 3 and 4 describe two large classes of applications.
In section 3, a first series of direct applications consists of models where the blocked 
edges themselves form the matter degrees of freedom. 
Section 3.1 deals with the statistics of Eulerian maps with bi- and tetra-valent
faces and with blocked edges, which can also be viewed as decorated quadrangulations whose
unoriented edges can be blocked in both directions independently, like maps
with possible one-way roads and/or fully blocking roadworks. The simpler case
of quadrangulations with fully blocked edges only is treated in section 3.2 where 
a detailed analysis of the critical singularities is presented. 
This allows in particular to solve the problem of quadrangulations with a spanning
tree. In appendix A, we provide for both solutions a simple combinatorial
interpretation in terms of maps with vertices of even valence and without matter.
Appendix B deals with the case of maximally blocked Eulerian maps, i.e. maps
with a maximal number of blocked edges but which remain connected in the 
sense that every vertex can still be reached via some path from the origin.

Section 4 is devoted to applications to maps with particles subject to
exclusion rules. We first show in section 4.1 how to formulate these problems
in terms of Eulerian maps with blocked edges by associating to any unconstrained configuration
of particles on a map a set of configurations of blocked edges on the same map with
possible negative weights, and whose total contribution vanishes as soon as the particle exclusion
rule is violated. This allows for the derivation in section 4.2 of recursive relations
for the generating functions of the models. Two particular examples are discussed:
Eulerian triangulations with hard particles in section 4.3 and, in section 4.4, the Ising
model on quadrangulations via a straightforward bijection with hard particles.
We gather a few concluding remarks in section 5.

\newsec{Coding maps with mobiles}

Eulerian maps are planar embedded graphs (considered up to continuous
deformation), which may be drawn in a single loop without lifting the pen. 
As each vertex is of even valence, the edges of these maps are naturally
oriented, with incoming and outgoing edges alternating around each vertex.
We may then imagine a walker moving along these oriented edges as on one-way
streets throughout the map.
We wish here to introduce the notion of Eulerian map with {\it blocked edges},
by which we simply mean that some of the edges of the map may be forbidden 
to the walker, like streets blocked by road works or so. The precise position of the 
blocked edges will be in general arbitrary up to one important restriction: 
we demand that every vertex of the map may be reached from a chosen origin vertex 
despite the blockings. 
As we shall see later, many statistical models on random maps may be described 
in terms of such Eulerian maps with blocked edges. This includes for instance 
spanning trees, hard particles, or even Ising spins.

\subsec{Eulerian maps with blocked edges}

Let us now come to more precise definitions of Eulerian maps with blocked edges.
An Eulerian map is equivalently defined as a planar map whose faces are colored
in black and white so that no two adjacent faces be of the same color. 
In other words, its dual map\foot{Recall that the dual of a planar map
is the planar map obtained by replacing each face by a vertex at its center and
each edge by a ``perpendicular" edge connecting the centers of the adjacent
faces. The vertices of the original map become the faces of the dual map.} 
is a bipartite planar map, namely
a map whose vertices may be partitioned into two subsets such that no two
vertices of the same subset are adjacent. The coloring of the 
faces of an Eulerian map induces a natural orientation of its edges by simply 
demanding for instance that this orientation be clockwise around black faces.
Here we shall consider so-called {\it pointed maps}, i.e. maps with a 
distinguished vertex, hereafter referred to as the {\it origin}.  
\fig{An example (a) of pointed Eulerian map with blocked edges (thick lines). 
The blockings satisfy the global connectivity constraint that every vertex can be 
reached from the origin by a path consisting of non-blocked edges taken in their 
canonical orientation (indicated by an arrow). On the dual bipartite map (b), 
the duals of blocked edges necessarily form a forest, i.e. a graph without cycles 
and, moreover, are such that each face can be reached from the external face 
(dual to the origin) without crossing the forest edges and by going clockwise 
(resp. anticlockwise) around black (resp. white) vertices.}{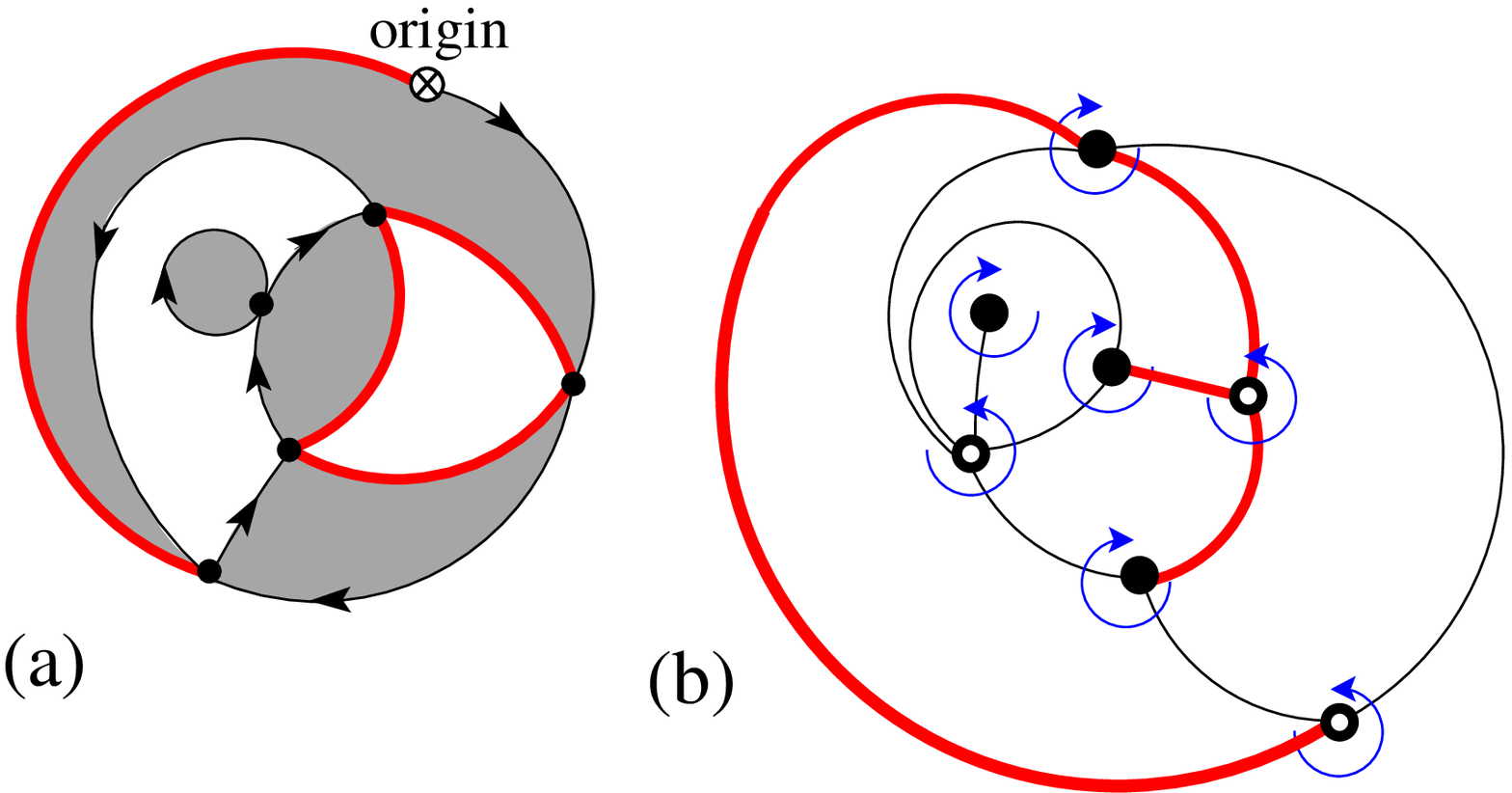}{12.cm}
\figlabel\eulermapwithbe
By map with blocked edges, we simply mean that we restrict the set of paths allowed
on the map by deciding that some of the edges of the map are ``blocked", while paths 
are only allowed to pass through un-blocked edges and must respect the orientation on 
these edges. 
As mentioned above, we furthermore impose a global {\it connectivity constraint} 
that each vertex be attainable from the origin by some allowed path (see figure 
\eulermapwithbe-(a) for an example).
{}From this requirement, we immediately see that, on the dual bipartite map, 
the set of edges dual to blocked edges (also referred to as blocked edges for
convenience) form a forest, i.e. contains no cycle. This forest is however
not arbitrary in general as it must be such that any face of the dual map can be reached from the 
origin face (dual of the origin vertex of the Eulerian map) without crossing blocked 
edges and by going clockwise (resp. anticlockwise) around black (resp.  white) vertices 
(see figure \eulermapwithbe-(b)).
The connectivity constraint is crucial for the construction below to be valid. 
In the following, we shall be led to introduce further restrictions on the blocked edge 
configurations, whose consequence on our construction will be transparent.

\subsec{Well-labeled mobile construction}

\fig{The rules for the construction of a well-labeled  mobile (see text). A flagged edge 
can be marked (thick line) or not (thin line) according to whether the dual edge
was blocked (i) or not (ii).}{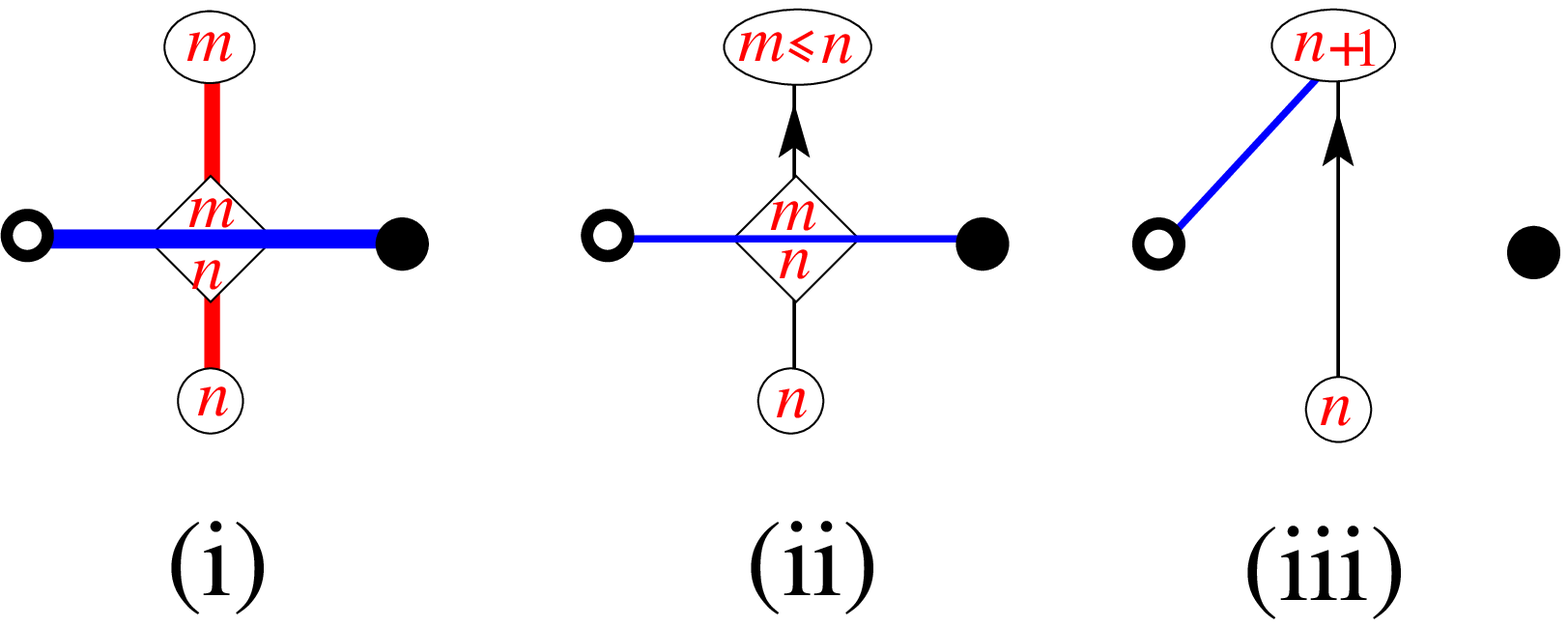}{12.cm}
\figlabel\mobconst
\fig{Construction of the well-labeled mobile associated with the Eulerian map
with blocked edges (a) of figure \eulermapwithbe. Applying the rules of figure \mobconst\ for 
each edge of the map (b) results in the desired mobile (c). On this mobile, flagged
edges can be marked (thick lines) or not.
}{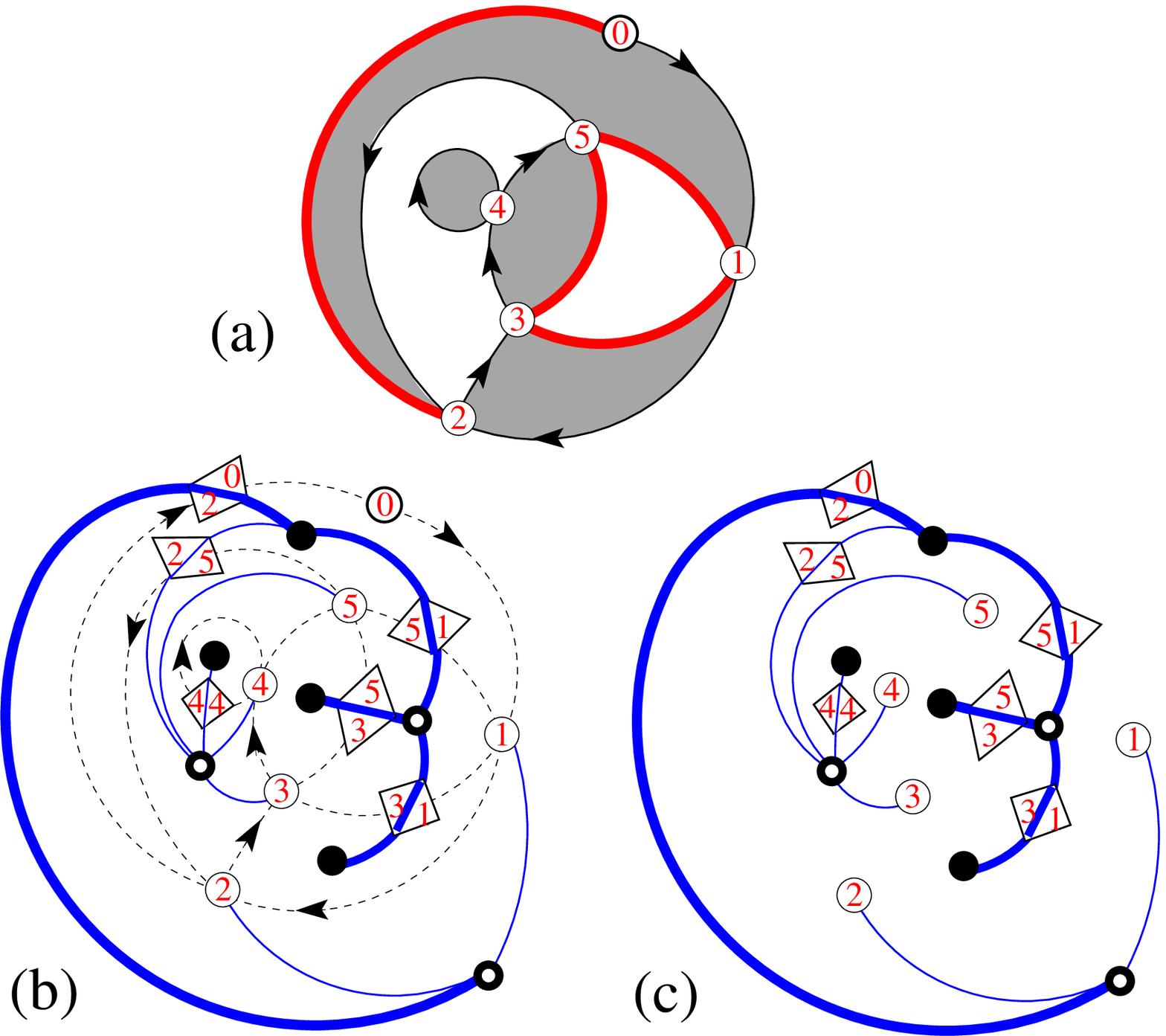}{12.cm}
\figlabel\mobileconst

The connectivity constraint above allows to assign to each vertex of the map
its (finite) {\it distance} from the origin defined as the minimal number of
steps of a path going from the origin to that vertex via un-blocked edges only 
and respecting the orientation.  Such a path will be referred to as a {\it geodesic} path.

We can encode the configuration of the map and its blocked edges via a so-called
{\it well-labeled mobile} as follows:

\item{1.} We first label all vertices of the map by their distance from the
origin.
\item{2.} We then add a white (resp. black) vertex at the center of each white
(resp. black) face.
\item{3.} We then consider each edge of the map and apply the following
construction:
\itemitem{(i)} If the selected edge is a blocked edge, we draw its dual edge by
connecting the black
and white vertices at the center of its adjacent faces and complete it with a
{\it flag} on
each side. We put on each flag the label of the closest vertex (see figure
\mobconst-(i)). 
We finally {\it mark} this dual edge. 
\item{}If the selected edge is not blocked, it points from a vertex at distance
$n$ from the
origin to a vertex at distance $m$, with necessarily $m\leq n+1$. 
\itemitem{(ii)} If $m\leq n$, we repeat the above construction and draw the dual
edge by connecting 
the black and white vertices at the center of its adjacent faces and complete it
with a flag on
each side, with again the same label as that of the closest vertex (see figure
\mobconst-(ii)). 
The dual edge is however un-marked in this case.
\itemitem{(iii)} If $m=n+1$, we draw an edge from the white vertex at the center
of the adjacent white
face to the vertex labeled $m=n+1$ (see figure \mobconst-(iii)).\par
\noindent The well-labeled mobile is simply defined as the graph obtained from
the collection of all the black 
and white vertices, all the added marked (case (i)) and un-marked (case (ii))
flagged edges with their 
labeled flags, and all the added edges from white vertices (case (iii)) together
with their labeled 
endpoints (see figure \mobileconst\ for an example). For convenience, we shall
refer 
to marked flagged edges, un-marked flagged
edges and edges between white vertices and labeled vertices as edges of type
(i), (ii) and (iii)
respectively. Note that, by construction, all labels on labeled vertices are
positive integers (since the origin vertex cannot be in the mobile),
while all labels on flags are non-negative integers.
The above construction extends that of Ref.\MOB\ that associates well-labeled
mobiles to 
Eulerian maps without blocked edges by use of (ii) and (iii) only. 

Let us now come to the characterization of the well-labeled mobiles obtained
from this construction.
The first and most important characterizing property of a well-labeled mobile is
that it is a tree, 
i.e. a connected graph without cycles. 
The proof is similar to that of Ref.\MOB. First we note that edges belonging to
a geodesic path yield only un-marked edges of type (iii) in the mobile, as the
distance increases by one at each step. 
This guarantees that geodesic paths on the map never cross the flagged edges of
the mobile, and 
may cross the mobile only at its labeled vertices. 
We then assume by contradiction the existence of a cycle on the mobile, and call
interior
the region delimited by this cycle that does not contain the origin. We then
consider the smallest 
label $n$ among the labels of the (labeled) vertices on the cycle and the labels
of the flags 
along the cycle that lie in the interior region.
If $n$ is attained for a labeled vertex on the cycle, we immediately deduce
from the construction (iii) the existence of a vertex at distance $n-1$ from the
origin 
in the interior. Any oriented geodesic path from the origin to that vertex must
intersect 
the cycle at a labeled vertex with label strictly smaller than $n$. This is a
contradiction.
If $n$ is not attained at a labeled vertex, but only at a flag,
we deduce from the construction (i) or (ii) the existence of vertex at distance
$n$ from the
origin in the interior, hence a geodesic path to that vertex intersects the
cycle at a
labeled vertex with label $\leq n$. This is again a contradiction.
We deduce that the mobile has no cycle, hence is a forest made of $c$ connected
components,
where $c$ is also the difference between its number of nodes and its number of
edges.
By construction, this number of edges is simply the number $E$ of edges of the
original 
map, while the number of nodes is equal to the number of faces $F$ of the
original map 
plus the number of labeled vertices.  
The latter is at most $V-1$ if we denote by $V$ the total number
of vertices of the original map, as the origin cannot belong to
the mobile (it remains isolated by construction hence is removed). As $F+V-E=2$
from the Euler relation for a planar map, we deduce that 
$c \leq F+(V-1)-E =1$, hence necessarily
$c=1$. The mobile is thus a tree that moreover contains all the vertices of
the original map but the origin.

\fig{Constraints on the labels around a white vertex of a mobile (see text).}{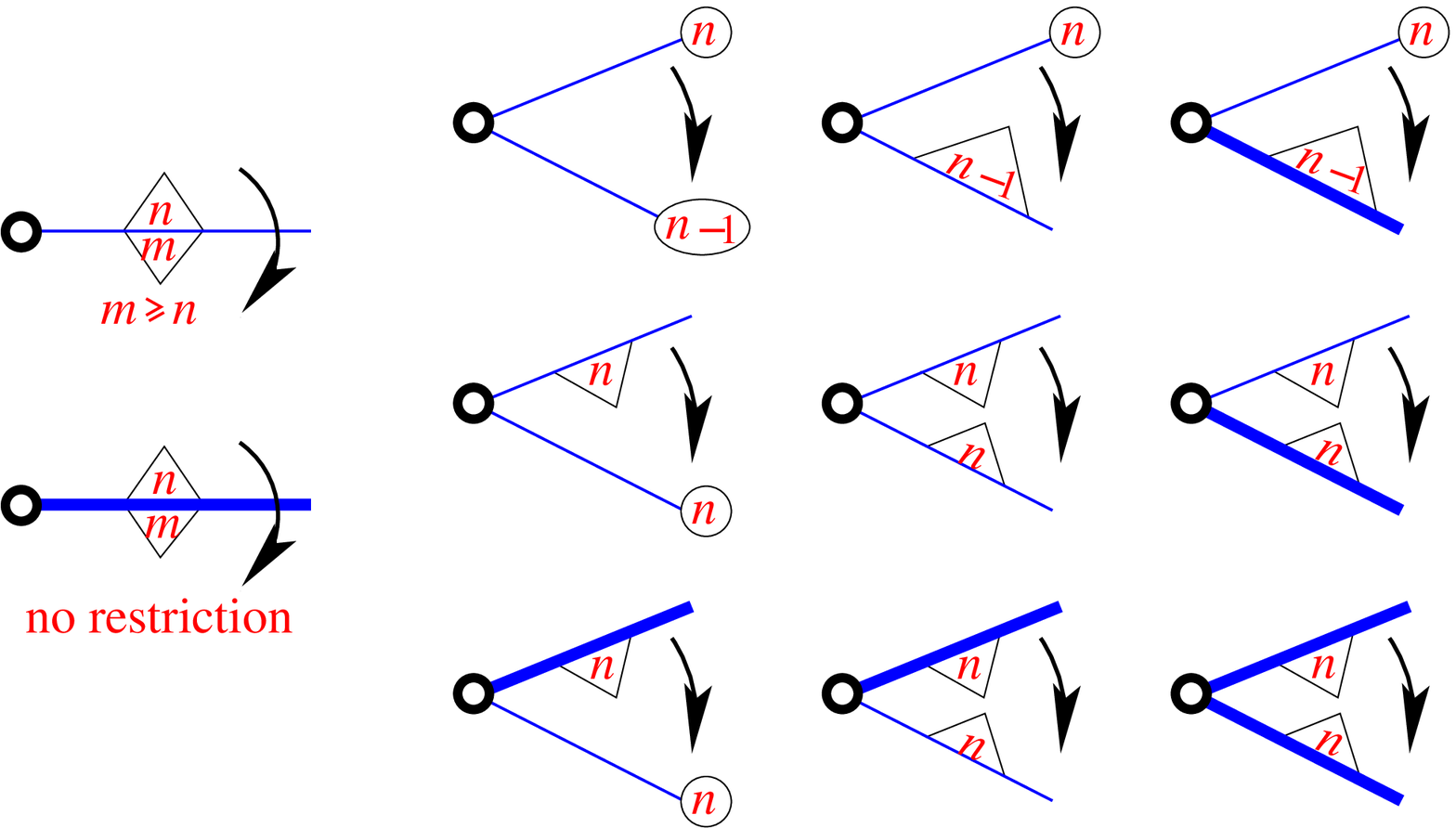}{12.cm}
\figlabel\white

The second characterizing property of well-labeled mobiles concerns the
environment of white
and black vertices, i.e.\  the configuration of labels around them. Recall that,
by {\it corner} at a vertex is meant any sector
delimited by two consecutive incident edges around that vertex.
In the case of a white vertex, we first note that, by construction,
the number of incident edges is nothing but the valence of the associated 
white face on the original map. We have the following rules for the cyclic
evolution of the 
labels around each white vertex (see figure \white):
\item{--} at the {\it crossing} of each un-marked flagged edge, the value of the
flag must increase
weakly clockwise. There is {\it no such restriction for marked flagged edges}.
\item{--} going through each {\it corner} clockwise around a white vertex, the
label decreases 
by one if the first edge of the corner is of type (iii), and remains constant if
the first
edge is of type (i) or (ii). \par
\fig{Constraints on the labels around a black vertex of a mobile (see text)}{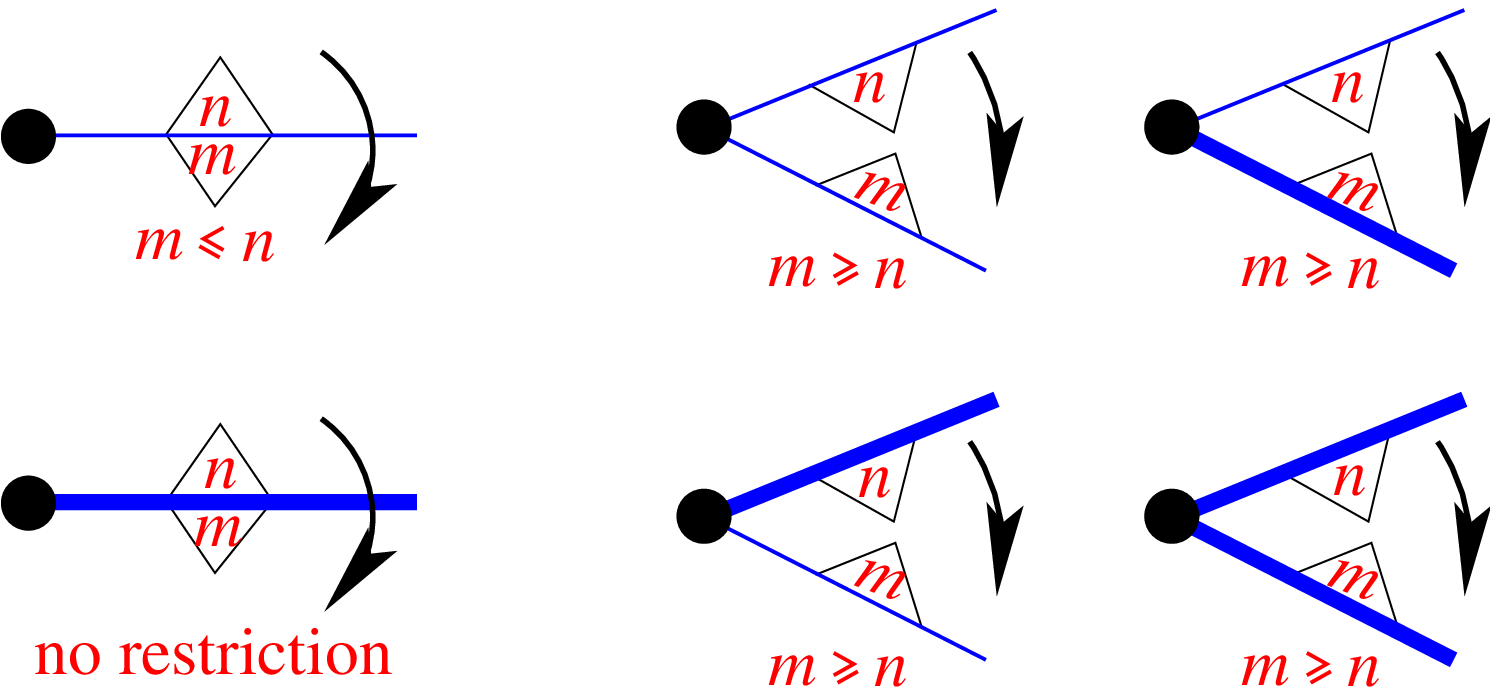}{10.cm}
\figlabel\black

\noindent As for black vertices, their incident edges are of type (i) and (ii)
only and we 
have the following 
rules for the cyclic evolution of the labels around each black vertex (see
figure \black):
\item{--} at the crossing of each un-marked flagged edge, the value of the
flag must decrease
weakly clockwise. There is no such restriction for marked flagged edges.
\item{--} going through each corner clockwise around a black vertex, the
label
increases weakly. \par
\fig{An alternative representation using spurious dangling edges of the constraints 
of figure \black\ for the evolution of labels at a corner around
a black vertex.}{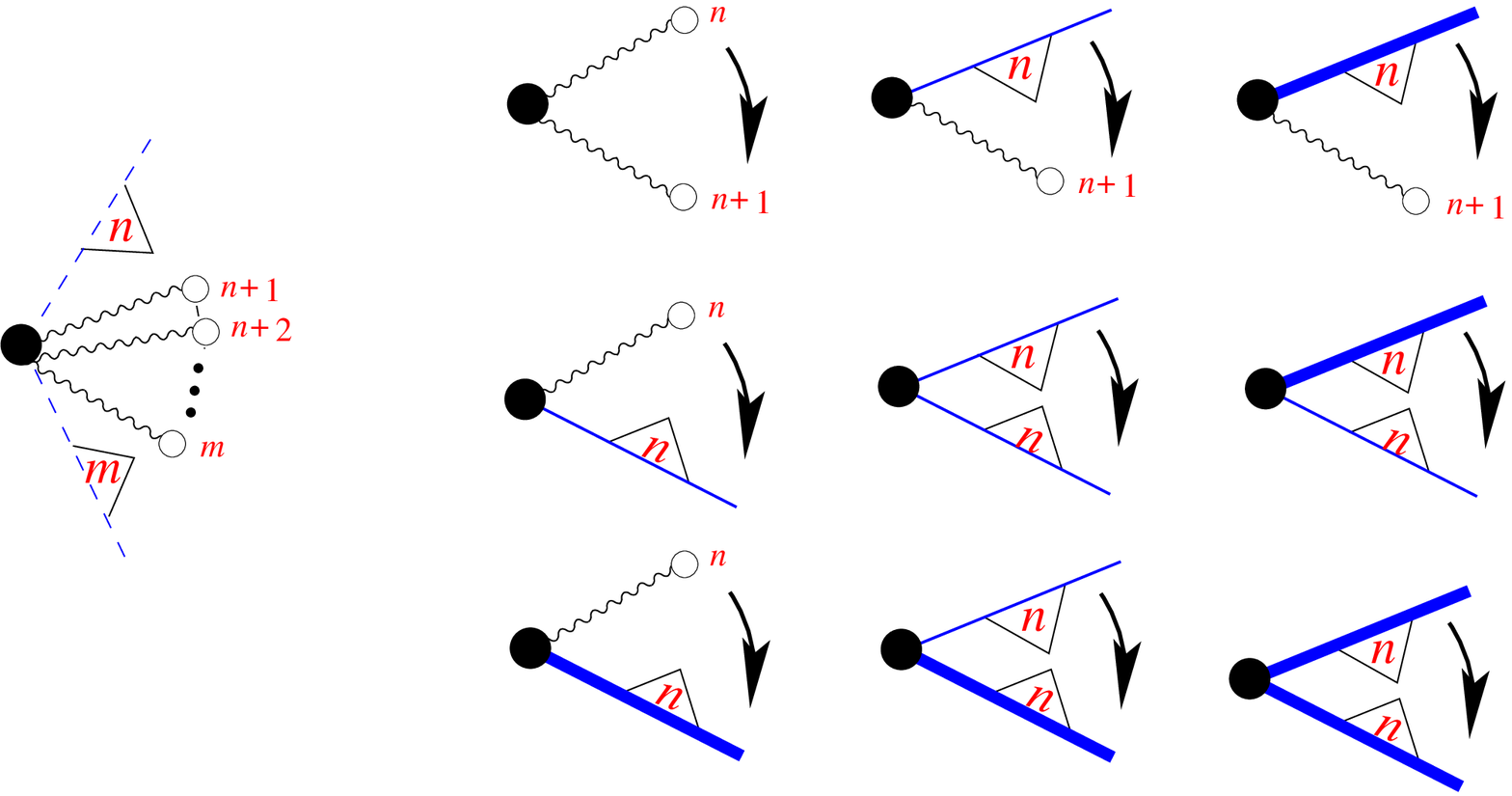}{12.cm}
\figlabel\blackbis
\noindent Note that the first (crossing) rules for white and black vertices are
actually redundant with each other, but we mention both of them to emphasize the
symmetry. Moreover, as displayed in figure \blackbis, the second (corner) rule
for black vertices can be made similar to 
that of white vertices by dressing
each corner around a black vertex with first label $n$ and second label $m\geq
n$ clockwise
by $m-n$ spurious dangling (wiggled) edges attached to spurious vertices labeled
$n+1, n+2, \ldots, m$. 
With this dressing, the black corner rule becomes:
\item{--} going through each corner clockwise around a black vertex, the
label increases
by one if the last edge of the corner is a dangling edge and remains constant if
the last
edge is of type (i) or (ii). \par
\noindent In this latter formulation, the total number of edges around a black
vertex, including the 
spurious ones, is nothing but the valence of the associated black face on the
original map.

The above properties for the evolution of labels, together with the fact that
these labels 
are positive integers on vertices and non-negative integers on flags, with at
least 
one vertex labeled $1$ or one flag labeled $0$, form a complete characterization
of the well-labeled 
mobiles. In the following we will be led to relax the above constraints on
labels. Mobiles for which we relax the constraint of existence 
of a label $1$ or $0$ but keep the positivity constraint will be referred to
as {\it positive mobiles}. Finally, mobiles for which we relax both constraints
will simply be referred to as {\it unrestricted mobiles}.

\subsec{Inverse construction}

The well-labeled mobiles constructed above form compact encodings of the
associated Eulerian maps with
blocked edges.  Indeed, to recover the original map and its blocked edges from
the corresponding mobile, we
simply apply the construction of Ref.\MOB\ as follows. First, by {\it labeled
corner} we simply mean a corner at a labeled vertex of the mobile, with the
corresponding label value. We then 
move clockwise {\it around the mobile} and record the cyclic sequence of
successive labels 
encountered at flags and/or labeled corners (spurious edges and vertices
displayed in figure \blackbis\ are not considered). The passage from a label to
the next is always via a
white or a black vertex, hence, from the rules obeyed by labels around those
vertices, we
immediately deduce that the sequence of labels decreases by one after each
labeled corner 
and increases weakly after each flag (in particular, it remains constant when
passing along 
a white vertex).  Otherwise stated, the sequence of labels around the mobile has
the ``ratchet'' 
property that it can have arbitrary increasing steps but decreasing steps of
{\it at most} $1$ 
which take place just after labeled corners. 

\fig{Inverse construction: starting from a well-labeled mobile (a), each 
labeled corner or flag is connected to its successor (b). The graph formed by these
chords (c) is the desired Eulerian map (d), where the blocked edges are obtained 
from the chords originating from marked flagged edges.}{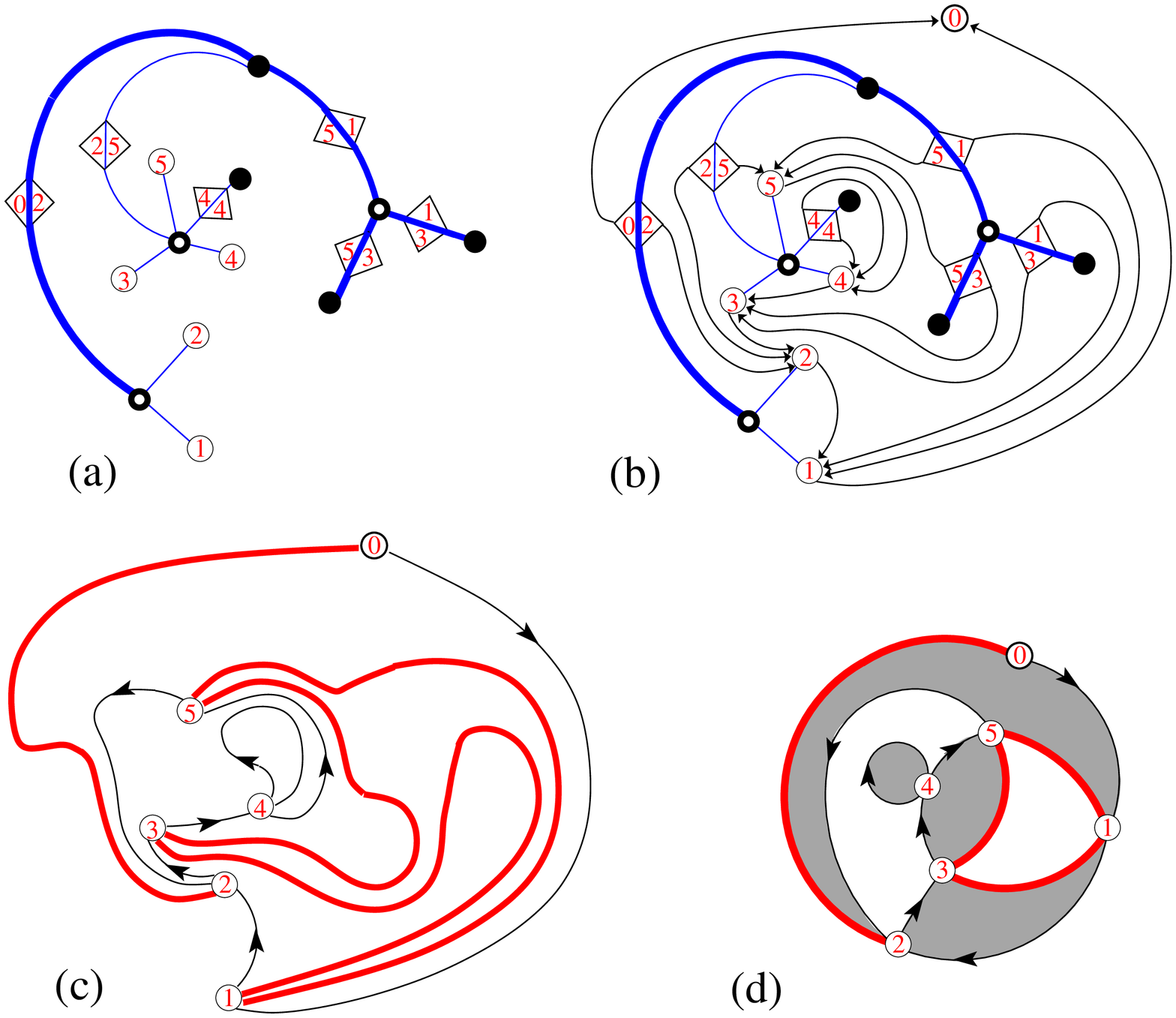}{12.cm}
\figlabel\backtoEuler
This property immediately ensures the existence of {\it successors} defined as
follows:
the successor of a labeled corner with label $n>1$ is the first labeled corner
encountered clockwise 
with label $n-1$. Similarly, the successor of a flag with label $n>0$ is the
first labeled corner 
encountered clockwise with label $n$. Note that the ratchet property of the
label sequence also
guarantees that if a corner or flag $l'$ lies between a corner or flag $l$ and
its successor $s(l)$, 
the successor $s(l')$ is necessarily attained at or before $s(l)$. We may
therefore connect each 
labeled corner with label $n>1$ and each flag with label $n>0$ to its successor
by a new edge, called a {\it chord}, in such 
a way that all these chords do not cross, while all corners labeled $1$ and all
flags labeled $0$ 
remain adjacent to the exterior face. We finally add an extra origin vertex
inside the exterior face
and connect all corners labeled $1$ and all flags labeled $0$ to that vertex by
non-crossing chords.
The original Eulerian map is made of the origin, the labeled vertices and  
the chords. Note that some chords (those starting from flags) have to be glued at
the level of the double flags so as to form 
the desired edges between the (origin or labeled) vertices of the map (see
figure \backtoEuler\ for an example).
Any edge of the Eulerian map is of one of the following three types (see figure
\backtoEuler): 
\item{(i)} either it is obtained by gluing two chords at the level of the double
flag of a marked edge.
We then make it a blocked edge.
\item{(ii)} or it is obtained by gluing two chords at the level of the double
flag of an un-marked edge.
This edge is then left un-blocked and we decide to orient it so that the black
endpoint of the flagged edge
lies on its right.
\item{(iii)} or it is obtained from a single chord connecting a labeled vertex to
its successor. This
edge is then left un-blocked and we decide to orient it backwards from the
successor. \par
\noindent The proof that the above construction and that of the previous section
are inverse of one another
is similar to that of Ref.\MOB\ in the absence of blocked edges and will not be
reproduced here.
The orientation of the edges simply reproduces that induced by the black and
white faces on the Eulerian map.
Finally, the labels on the vertices simply record the distance of these vertices
form the origin 
by oriented paths avoiding blocked edges. Indeed, form our choice of orientation
above,
the sequence of labels along any oriented path from the origin on the map either
decreases 
weakly (case (ii) above) or increases by one (case (iii) above) at each step.
The label of a vertex
is therefore necessarily smaller than its distance from the origin. In other
words, the distance from
the origin to a vertex is larger or equal to its label, and it must be exactly
equal as the sequence
of consecutive successors of that vertex provides, when considered backwards, an
oriented path 
from the origin to that vertex whose length is exactly equal to its label.

\subsec{Generating series}

The main interest of the above coding by mobiles concerns the {\it enumeration}
of families of 
Eulerian maps with blocked edges. Indeed, such enumeration simply
amounts to counting appropriate mobiles and the underlying tree structure 
of these mobiles translates into simple recursive equations for the associated
generating 
functions. 

In the following, we shall consider only mobiles with at least one labeled
vertex.
In other words, we exclude the degenerate case of mobiles made only
of flagged edges, necessarily all with equal labels. Going back to maps,
this amounts to excluding the quite trivial maps made of a single vertex only,
connected 
by a configuration of nested ``petals''. The enumeration of mobiles is then made
easier by considering {\it rooted} mobiles, 
i.e. mobiles with a distinguished corner around a labeled vertex.  We may then
view the mobile as hanging from an extra root edge added 
in this corner. 

Let us now consider the quite general problem of enumerating
maps with {\it prescribed face valences}
and with a fixed number of blocked edges. At this stage, we make for simplicity
no restriction on the configuration of blocked edges around a face other than
the global connectivity constraint. Some extra restrictions will be introduced
in some applications below, whose incorporation will lead to straightforward and
harmless modifications. As usual, it is simpler to work within grand 
canonical ensembles of maps in which the numbers of faces of given valences may vary 
and are controlled by
weight factors. More precisely, we consider the set of maps with faces of
arbitrary 
valences up to some upper bound $p$, counted with weight $g_k$ per $k$-valent
white face, 
$\tilde g_k$ per $k$-valent black one ($k=1,2,\ldots,p$). We also let the
number of blocked edges vary, with a weight $y$ per blocked edge.
This amounts to consider mobiles with white and black vertices of valence $k\le
p$ (including the spurious dangling edges for black vertices), weighted
respectively by 
$g_k$ and $\tilde g_k$, and with an arbitrary number of marked edges weighted
by $y$. 

We denote by $R_n$ the generating function of such weighted mobiles rooted at a
corner labeled $n$. At this stage, we need not specify whether the labels obey
positivity conditions (case of positive mobiles) or not (case of unrestricted
mobiles).
Indeed, the following reasoning holds in both cases, the difference coming only
from boundary conditions that we shall discuss later. We shall not consider the
case of well-labeled mobile (rooted at a corner labeled $n$) as their generating 
function $G_n$ is directly expressed in terms of that for positive mobiles 
via $G_n=R_n-R_{n-1}$. 
\fig{Recursive decomposition of a mobile rooted at a labeled corner. When the root vertex
has valence $k$, the decomposition yields $k$ mobiles rooted at a univalent labeled vertex.
In terms of generating function, this translates into eq.(2.1).
}{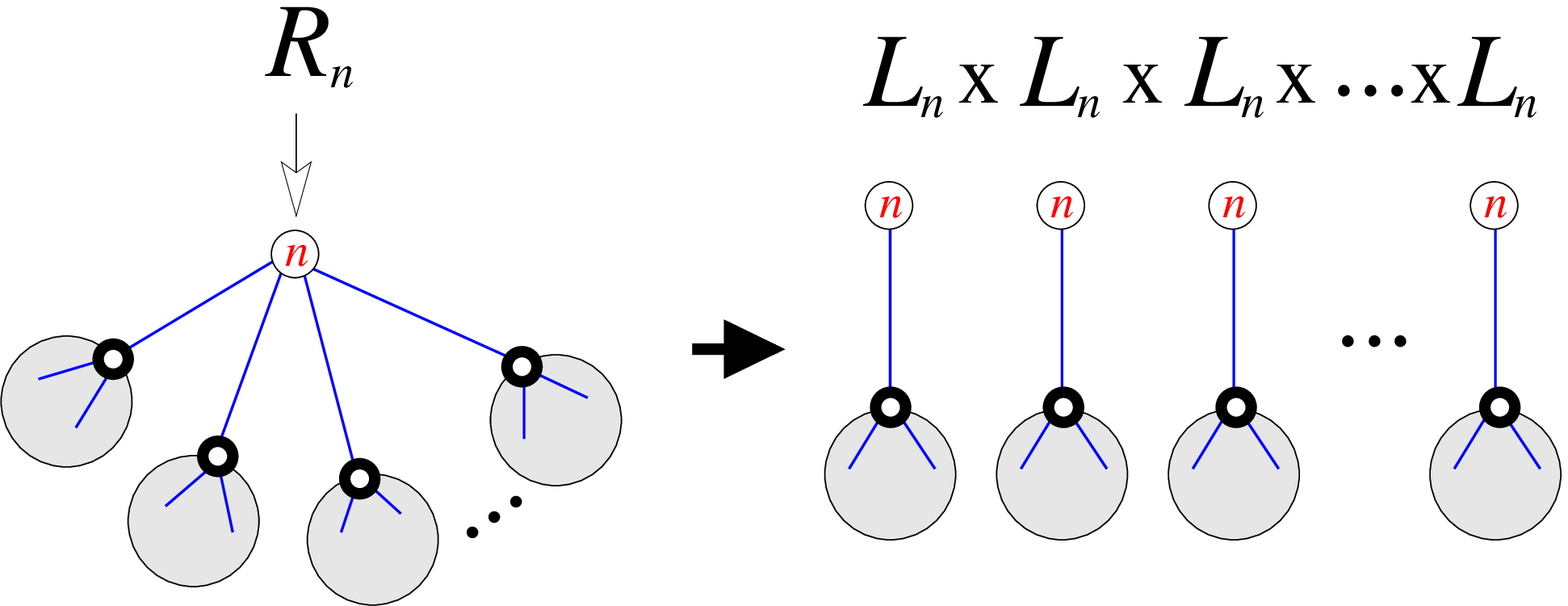}{10.cm}
\figlabel\recursion
Let us derive the set of recursion relations for $R_n$ inherited from the mobile
structure.
We have the first relation 
\eqn\RtoL{R_n=\sum_{k=0}^{\infty} (L_n)^k={1\over 1-L_n}\ ,}
where $L_n$ is the generating function of mobiles rooted at a {\it univalent}
vertex
labeled $n$. This is readily seen by decomposing the root vertex according to
its arbitrary valence $k$ (see figure \recursion). 
\fig{The three possible types of descending subtrees at a white vertex of a mobile according
to the type of their root edge, together with the associated generating functions.}{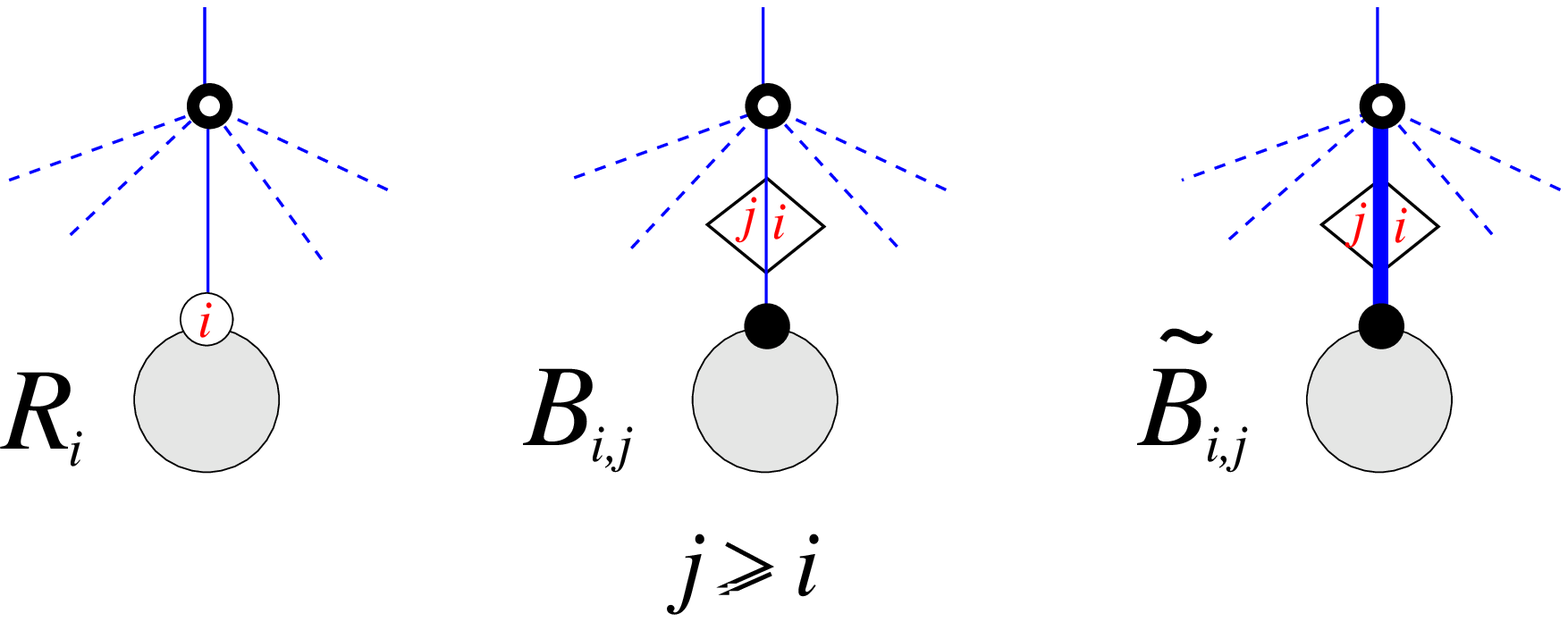}{10.cm}
\figlabel\whitesubtree
\fig{The three possible types of descending subtrees at a black vertex of a mobile according
to the type of their root edge, together with the associated generating functions. It proves convenient 
to use the representation of fig.\blackbis\ of mobiles with spurious dangling edges 
(with a trivial generating function $1$) which emphasizes the similarity between black
and white vertices.}{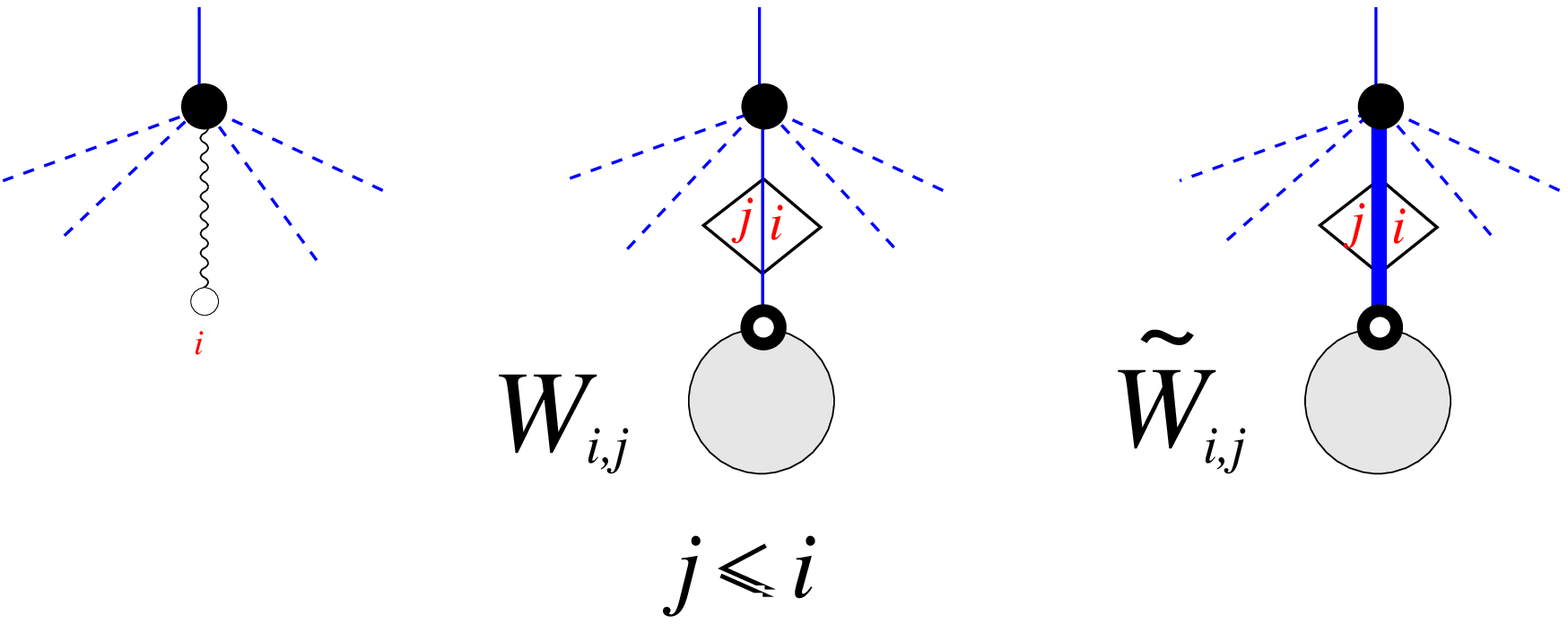}{10.cm}
\figlabel\blacksubtree

The generating function $L_n$ is then obtained by inspecting
the environment of the white vertex attached to the univalent root vertex. Its
descendent subtrees come into three species (see figure \whitesubtree):
\item{(i)} sub-mobiles rooted at a corner around a labeled vertex, with
generating function $R_i$ if the root label is $i$, 
\item{(ii)} sub-mobiles rooted at some un-marked 
flagged edge incident to a black vertex, with generating function $B_{i,j}$ if
$i$ (resp. $j$)
is the label on the right (resp. left), viewed with the white vertex on top,
\item{(iii)} sub-mobiles rooted at some marked flagged edge incident to a black
vertex, with generating 
function $\tilde B_{i,j}$ if again $i$ (resp. $j$) is the label on the right
(resp. left). \par
\noindent To follow the evolution of labels clockwise around the inspected
white vertex, we introduce a formal (transfer) operator ${\bf Q}_\circ$ acting
on basis 
vectors $|i\rangle$
indexed by the encountered labels. The action of this operator reads:
\eqn\Qcirc{{\bf Q}_\circ|i\rangle = R_i|i-1\rangle\, +\sum_{j\geq i}
B_{i,j}|j\rangle 
+\sum_{j} \tilde B_{i,j}|j\rangle \ ,}
where the ranges of summation incorporate the white crossing rules (obeyed by
labels at the crossing of flagged
edges) displayed in fig.\white. 
Returning to the computation of $L_n$, we note that, as the first label encountered 
clockwise is $n-1$, the succession of labels around the
white vertex must lead from $n-1$ back to $n$ in $k-1$ steps if the white vertex
is $k$-valent.
This leads to
\eqn\LtoB{L_n=\sum_k g_k \langle n | {\bf Q}_\circ^{k-1} | n-1 \rangle\ ,}
where $\langle i|$ denotes a vector of the dual basis defined by $\langle
i|j\rangle=\delta_{i,j}$.

Likewise, inspecting the environment of a black vertex incident to a flagged
root edge,
we get 
\eqn\BtoW{B_{i,j}=\sum_k \tilde{g}_k \langle j | {\bf Q}_\bullet^{k-1} | i
\rangle\ ,\ j\geq i }
and  
\eqn\BtildetoW{\tilde B_{i,j}=y\,\sum_k \tilde{g}_k \langle j | {\bf
Q}_\bullet^{k-1} | i \rangle\ ,}
where the transfer operator ${\bf Q}_\bullet$ around a black vertex now reads
\eqn\Qbullet{{\bf Q}_\bullet|i\rangle = |i+1\rangle 
+ \sum_{j\leq i} W_{i,j}|j\rangle
+ \sum_{j} \tilde W_{i,j}|j\rangle.
}
Here the three terms account for the three species of sub-mobiles one may
encounter
around a black vertex (see figure \blacksubtree):
\item{(i)} dangling ends with weight $1$, 
\item{(ii)} sub-mobiles rooted at some un-marked 
flagged edge incident to a white vertex, with generating function $W_{i,j}$ if
$i$ (resp. $j$)
is the label on the right (resp. left), viewed with the black vertex on top,
\item{(iii)} sub-mobiles rooted at some marked flagged edge incident to a white
vertex, with generating 
function $\tilde W_{i,j}$ if again $i$ (resp. $j$) is the label on the right
(resp. left). \par
\noindent The ranges of summation in eq.\Qbullet\ now incorporate the black
crossing rules 
(obeyed by labels at the crossing of flagged edges) displayed in fig.\black.
The generating functions $W_{i,j}$ and $\tilde W_{i,j}$ are finally obtained
from:
\eqn\WtoB{W_{i,j}=\sum_k {g}_k \langle j | {\bf Q}_\circ^{k-1} | i \rangle\ ,\
j\leq i}
and
\eqn\BtildetoW{\tilde W_{i,j}=y\,\sum_k {g}_k \langle j | {\bf Q}_\circ^{k-1} |
i \rangle\ .}

Let us now discuss the distinction between positive and unrestricted mobiles. 
In the case of positive mobiles, the equations \RtoL-\BtildetoW\ above are only
valid 
for the allowed range of labels, i.e. $n>0$ for $R_n$ and $L_n$ and $i,j\geq 0$
for $B_{i,j}$, $\tilde B_{i,j}$, $W_{i,j}$ and $\tilde W_{i,j}$. 
Moreover, in eqs. \Qcirc\ and \Qbullet, we must implicitly assume that $R_i=0$
for $i\leq 0$ 
and $B_{i,j}=\tilde B_{i,j}=0$ as well as $W_{i,j}=\tilde W_{i,j}=0$ if $i<0$ or
$j<0$.
With these conditions, the set of equations \RtoL-\BtildetoW\ determine
completely all generating functions 
as formal power series of the $g_k$'s and the $\tilde g_k$'s. 

The case of unrestricted labels is much simpler as \RtoL-\BtildetoW\ are valid
for all
values of the indices and moreover, enjoy a manifest translation invariance
property. 
More precisely, we may write
\eqn\limitgen{\eqalign{
R&\equiv R_n\cr
B_\ell & \equiv  B_{n,n+\ell}\ \ \ell\geq 0 \cr
\tilde B_\ell & \equiv  \tilde B_{n,n+\ell }\cr
W_\ell & \equiv  W_{n,n+\ell }\ \ \ell\leq 0 \cr
\tilde W_\ell & \equiv \tilde W_{n,n+\ell }\cr
}}
as all the above quantities on the right hand side are independent of $n$.

Going back to maps, one can easily check that $R$ is the generating function for
Eulerian 
maps with blocked edges, with an origin vertex and with a {\it distinguished
edge} which is not blocked 
and moreover points from a vertex at distance $m$ from the origin to one at
distance $m+1$, for some $m$. 
Indeed, we have clearly a bijection between the set of well-labeled mobiles
rooted at a corner with
{\it arbitrary} label, and the set of unrestricted mobiles rooted at a corner
with {\it prescribed} label (say $0$),
whose generating function is $R$.
The bijection simply consists in shifting all labels of the well-labeled mobile
by a constant value, such that the root label becomes $0$. Finally, the corners
labeled $m+1$ on a well-labeled mobile 
are in correspondence with non-blocked edges of type $m\to m+1$ on the
corresponding map.

One may as well distinguish an edge on the map which is not blocked, but now points from
a label at distance $m$ from the origin to a label at distance $m-\ell$
($\ell\geq 0$)
for some arbitrary $m$. The corresponding generating function is easily seen to
be $B_\ell W_{-\ell}$. 
Finally, one may also distinguish an edge which is blocked and connects 
a vertex at distance $m$ from the origin to a label at distance $m-\ell$ 
for some arbitrary $m$. The corresponding generating function reads
$(1/y)(\tilde B_\ell \tilde W_{-\ell}+\tilde B_{-\ell} \tilde W_{\ell})$ for
$\ell \neq  0$
and $(1/y)\tilde B_0 \tilde W_{0}$ for $\ell=0$. 
To summarize, the generating function of pointed Eulerian maps with blocked
edges and
with a {\it distinguished non-blocked edge} is
\eqn\nbe{R+\sum_{\ell \geq 0} B_\ell W_{-\ell}}
while that of pointed Eulerian maps with blocked edges and
with a {\it distinguished blocked edge} is
\eqn\be{{1\over y}\sum_{\ell} \tilde B_\ell \tilde W_{-\ell} \ .}

The functions $R$, $B_\ell$, $W_\ell$, ${\tilde B}_\ell$ and ${\tilde W}_\ell$ 
are the translation invariant solutions of the equations \RtoL-\BtildetoW,  
which can be rewritten in a slightly more compact form by introducing a formal
variable $z$ 
and defining the Laurent series $Q_\circ(z)$ and $Q_\bullet(z)$ via 
\eqn\Qscalar{\eqalign{
Q_\circ(z) &= {R\over z}\, +\sum_{\ell\geq 0} B_{\ell}\, z^\ell
+\sum_{\ell} \tilde B_{\ell}\, z^\ell \cr
Q_\bullet(z) &=  z\, +\sum_{\ell\leq 0} W_{\ell}\, z^{\ell}
+\sum_{\ell} \tilde W_{\ell}\, z^{\ell}\cr
}}
that are translation invariant versions of eqs.\Qcirc\ and \Qbullet.
The recursive equations read:
\eqn\recursimp{\eqalign
{ R &= {1\over 1-L}\ ,\quad L=\sum_k g_k Q_\circ (z)^{k-1}|_{z^1}\cr
B_\ell &= \sum_k \tilde g_k Q_\bullet(z)^{k-1}|_{z^\ell}, \quad \ell\geq 0\cr
\tilde B_\ell &= y\ \sum_k \tilde g_k Q_\bullet(z)^{k-1}|_{z^\ell}\cr
W_\ell &= \sum_k g_k Q_\circ(z)^{k-1}|_{z^\ell}, \quad \ell \leq 0\cr
\tilde W_\ell &= y\ \sum_k g_k Q_\circ(z)^{k-1}|_{z^\ell}\cr
}}
where $|_{z^\ell}$ stands for the extraction of the coefficient of $z^\ell$ in
the corresponding Laurent series.

Eliminating $L$, all $B$'s and $W$'s, we arrive at the closed system of
equations 
\eqn\system{\encadremath{\eqalign{
Q_\circ(z) &= {R\over z} + \sum_k \tilde g_k \left[ Q_\bullet^{k-1}(z) \right]_+
 + y\, 
\sum_k \tilde g_k Q_\bullet^{k-1}(z)\cr
Q_\bullet(z) &= z + \sum_k g_k \left[ Q_\circ^{k-1}(z) \right]_- + y\,
\sum_k g_k Q_\circ^{k-1}(z)\cr
R& =1+\sum_k g_k \left[{Q_\circ^{k-1}(z)\over z}\right]_0  R \ ,\cr
}}}
where the notations $[\cdot]_+$, $[\cdot]_-$ and $[\cdot]_0$ stand respectively
for the extraction of the non-negative 
power part, the non-positive power part, and the constant term in the Laurent
series. In practice, this system may be solved in two steps: the first two lines
determine all coefficients of the Laurent series as functions of $R$ and 
formal power series of the $g$'s and $\tilde g$'s. Substituting these in
the third line yields a self-consistent equation for $R$, which itself is solved
order by order in the $g$'s and $\tilde g$'s.

When $y=0$ we recover the usual equations for Eulerian maps (without blocked
edges), as obtained from the mobile construction of Ref.\MOB\ or by the so-called
orthogonal polynomial solution of the two-matrix model. It is interesting to note
that the restriction of Laurent series to their positive or negative parts simply reflects the
constraints on flag labels for un-marked edges (crossing rules). For marked
edges (corresponding to the $y$ term) there is no such restriction. 

Finally, we have the {\it duality} property that, whenever $g_k = \tilde g_k$ for
all $k$, $Q_\circ(z)$ is equal to $Q_\bullet({R \over z})$, while $R$ can be alternatively
recovered via:
\eqn\alterR{R = 1 + \sum_k g_k \left[ Q_\bullet(z)^{k-1} z \right]_0 \ . }

\newsec{Direct applications}

In this section, we show for illustration how to use the above formalism to study 
the statistics of blocked edges in a number of simple cases in which the blocked edges
themselves form the matter degrees of freedom on the maps. 

\subsec{Blocked edges on quadrangulations viewed as Eulerian maps}

\fig{An example of pointed Eulerian map (left) with blocked edges (thick lines) satisfying
the global connectivity constraint, where all black faces are bivalent and all white 
faces tetravalent. Each black face can be squeezed into an edge, resulting in a quadrangulation
(right) whose edges can be either non-blocked (thick lines), blocked in one direction
only (thick lines with one-way signs) or blocked in both directions (thin lines with
roadwork signs).}{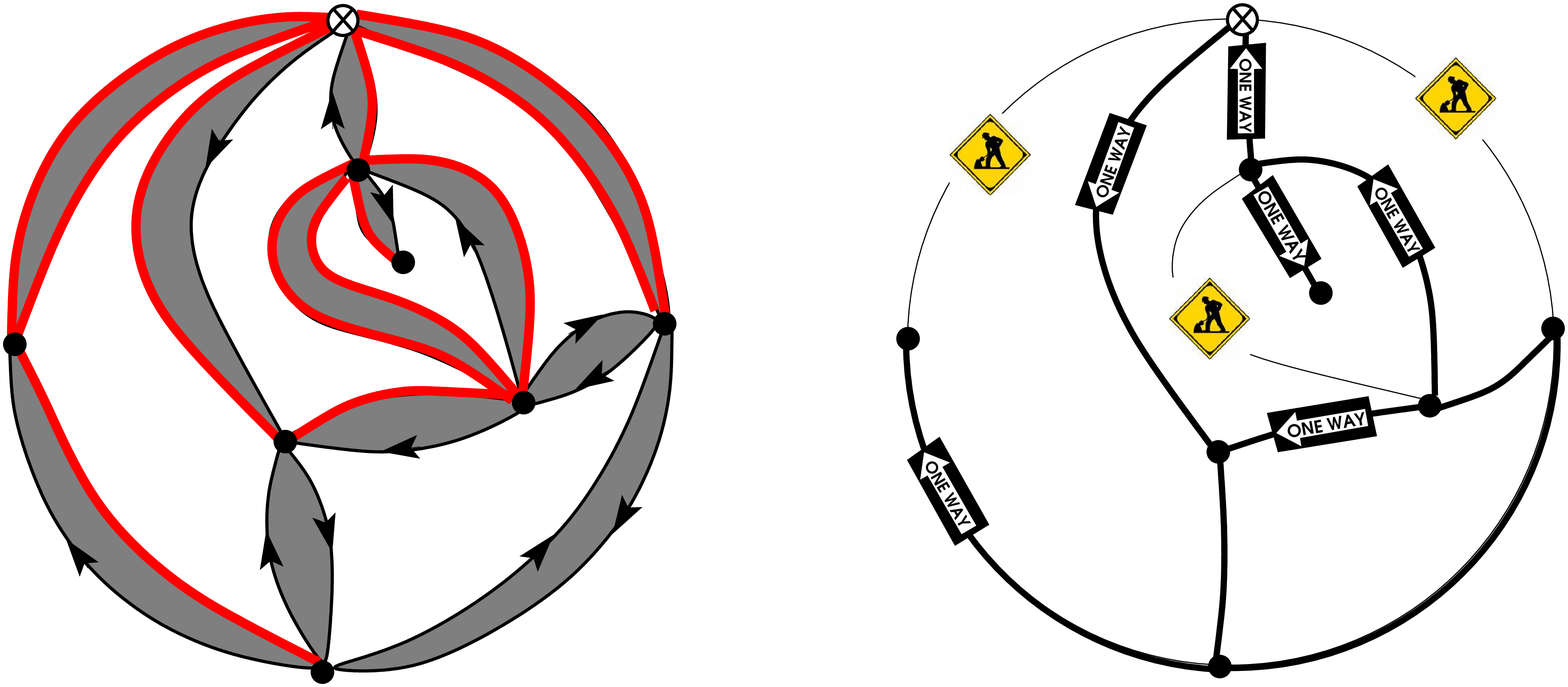}{13.cm}
\figlabel\quadran

As a first direct application of the results of section 2, let us consider the simple case
of quadrangulations viewed as Eulerian maps as follows: consider Eulerian maps where all
white faces are tetravalent and all black faces are bivalent. We may squeeze all the black 
faces into single edges, which results in a (non necessarily Eulerian) map with tetravalent
faces only, i.e. a quadrangulation. After this transformation, every edge replaces two edges of
the original map, with opposite orientations. Hence, in the absence of blockings, each edge
of the quadrangulation can now be taken in both directions. The introduction of blockings 
now corresponds to forbidding 
either one direction (one way road) or both, keeping the global connectivity constraint that 
every vertex be reachable from the origin of the map (see figure \quadran\ for an example). 
Taking  eq.\system\ with $g_k= g\, \delta_{k,4}$ and 
$\tilde g_k=\delta_{k,2}$ (as there are exactly twice as many black bivalent faces as white
tetravalent ones, we need not introduce an extra weight factor for the black faces), we get 
the system:
\eqn\firstappli{\eqalign{
Q_\circ(z) &= {R\over z} + \left[ Q_\bullet(z) \right]_+
 + y\, 
Q_\bullet(z)\cr
Q_\bullet(z) &= z + g\, \left[ Q_\circ^{3}(z) \right]_- +
g y\,  Q_\circ^{3}(z)\cr
R& =1+ g \left[{Q_\circ^{3}(z)\over z}\right]_0  R \ .\cr
}}
{}First we note that, by parity, both $Q_\circ(z)$ and $Q_\bullet(z)$ contain only odd
powers of $z$. From the second line, we then immediately deduce that
\eqn\qbulpos{\left[ Q_\bullet(z) \right]_+=z +  g y\, \left[ Q_\circ^{3}(z) \right]_+}
and, upon substituting in the first line
\eqn\eqforq{Q_\circ(z)= {R\over z}+(1+y)z+ g y (1+y) Q_\circ(z)^3 \ .}
Upon setting $x=g y (1+y)(R/z+(1+y)z)^2$ and $Q_\circ(z)=(R/z+(1+y)z)q(x)$,  
this equation becomes 
\eqn\eqq{q(x)=1+x q(x)^3} 
whose formal solution is
\eqn\solquadr{q(x)= \sum_{n\geq 0} {(3n)!\over n!(2n+1)!}x^n \ .}
Writing the last line of eq.\firstappli\ as
\eqn\lastline{\eqalign{R&=1+g \left[{1\over g y (1+y)}\left(Q_\circ(z)-{R\over z}-(1+y)z\right)
{1\over z}\right]_0  R\cr
&= 1+ \sum_{n\geq 1} {(3n)!\over n!(2n+1)!} g^n
\left(y(1+y)\right)^{n-1}\left[\left( {R\over z}+(1+y)z\right)^{2n+1}{1\over z}\right]_0 R 
\ ,\cr}}
we end up with a single equation for the generating function $R$:
\eqn\eqforR{R=1+ \sum_{n\geq 1} {(3n)!\over n!n!(n+1)!} g^n y^{n-1}
(1+y)^{2n}R^{n+1} \ .}
Rather than directly analyzing this formula, we will derive and study in the next section
a very similar and slightly simpler formula (see eq.(3.10) below). Their connection is 
discussed in appendix A below.

\subsec{Blocked edges on quadrangulations and spanning trees on the dual}
\fig{A quadrangulation (left) with edges (thick lines) blocked in both directions, 
and satisfying the global connectivity constraint (irrespectively of the choice of the origin
vertex). On the dual tetravalent map (right), duals of blocked edges (thick lines) 
form a forest.}{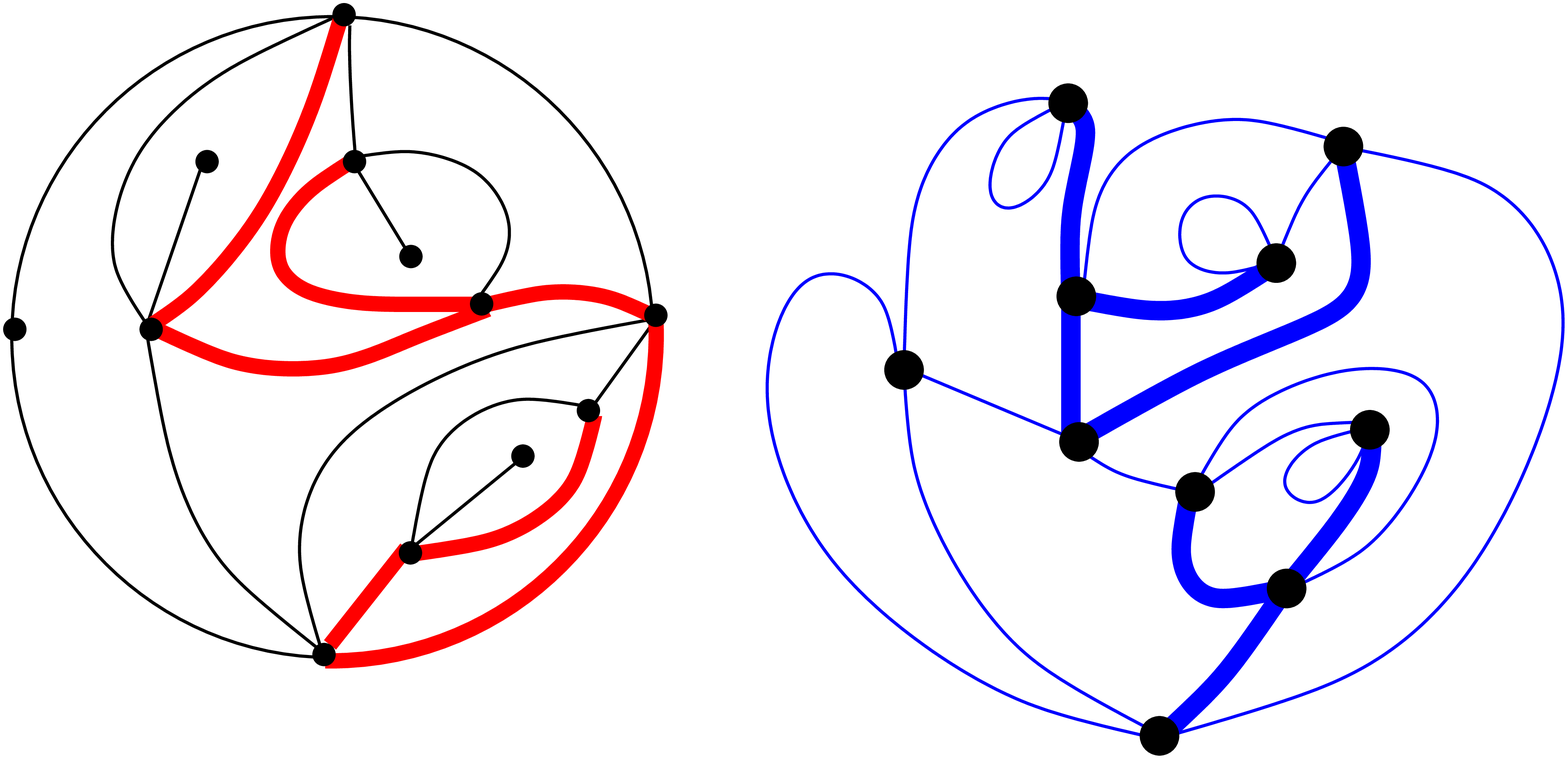}{13.cm}
\figlabel\quadranbis
As there is no canonical orientation of the edges in quadrangulations, a more
natural definition of blocked edges consists in either allowing for {\it both directions} 
on a given edge or forbidding both simultaneously.  
In the Eulerian map equivalent formulation of previous section, this simply means that
blockings
always take place for both directions simultaneously. In this case, for a fixed map, 
the configurations of blocked edges satisfying the global connectivity constraint turn out
to be independent of the choice of the origin. Rephrased in the dual language, this constraint
simply amounts to demand that the duals
of blocked edges make no cycles, i.e. form a forest (see figure \quadranbis\ for an example)
with no further restriction. 
The problem therefore amounts to enumerating {\it forests on planar maps} with tetravalent vertices,
with a weight $g$ per vertex and, say $y$ per edge of tree in the forest. 

Imposing that both direction be blocked simultaneously is achieved in the mobile language 
by requiring that every (bivalent) black vertex on the mobile has either $0$ or $2$ 
incident marked edges. This constraint results in a slight modification of eq.\firstappli,
now replaced by: 
\eqn\modifappli{\eqalign{
Q_\circ(z) &= {R\over z} + \left[ {\bar Q}_\bullet(z) \right]_+
 + y\, 
\tilde Q_\bullet(z)\cr
{\bar Q}_\bullet(z) &= z + g\, \left[ Q_\circ^{3}(z) \right]_- \cr
\tilde Q_\bullet(z) &= g \,  Q_\circ^{3}(z) \cr
R& =1+ g \left[{Q_\circ^{3}(z)\over z}\right]_0  R \ ,\cr
}}
which amounts to a splitting of the operator $Q_\bullet(z)$ of eq.\firstappli\ into
${\bar Q}_\bullet(z)$ and $\tilde Q_\bullet(z)$ respectively corresponding to sub-mobiles
with an un-marked and a marked root edge. 
The first three equations now reduce to 
\eqn\modifeqforq{Q_\circ(z)= {R\over z}+z+ g y\, Q_\circ(z)^3}
from which we deduce the equivalent of eq.\eqforR\ for the generating function $R$:
\eqn\eqforest{R=1+ \sum_{n\geq 1} {(3n)!\over n!n!(n+1)!} g^n y^{n-1}
R^{n+1}\ .}
For instance, the first few terms in the expansion of $R$ as power series of $g$ read
\eqn\expanRforest{\eqalign{R& = 1\!+\!3g\!+\!6g^2(3\!+\!5y)\!+\!15g^3(9\!+
\!30y\!+\!28y^2)\!+\!18g^4(63\!+\!315y\!\cr &+\!570y^2\!+\!385y^3)\!+\!126g^{5}
(81\!+\!540y\!+\!1440y^2\!+\!1855y^3\!+\!1001y^4)\cr &+\!36 g^{6}(2673\!+
\!22275y\!+\!78300y^2\!+\!146970y^3\!+\!149884 y^4\!+\!68068y^{5}) 
\!+\!{\cal O}(g^{7})\ ,\cr}}
from which we read off the generating functions for forests on tetravalent maps 
(with distinguished face and unoccupied edge) with up to $6$ vertices. 
In appendix A, we show how to reinterpret eq. \eqforest\ in the language of
even-valent maps without matter. 
Note that equation \eqforest\ is equivalent to eq.\eqforR\ upon substituting $g\to g(1+y)^2$.
The origin of this extra factor is also explained combinatorially in appendix A.

Let us now investigate the analytical behavior of $R$ as a function of $g$ for a 
fixed finite positive value of $y$ (note that each term of the series expansion
of $R$ in powers of $g$ is a polynomial in $y$). Introducing the function
\eqn\Fx{F(u)\equiv \sum_{n\geq 1} {(3n)!\over n!n!(n+1)!}u^n \ ,}
we can rewrite eq.\eqforest\ in a parametric form as
\eqn\parametric{\eqalign{g& ={u \left( y-F(u) \right) \over y^2} \cr
R&= {y\over y-F(u)}\ ,\cr}}
where we have introduced for convenience a parameter $u=g\,y\,R$. In this parametrization, 
$u$ must lie between $0$ and $1/27$, which is the radius of convergence of $F(u)$. However, 
the first equation determines implicitly $u$ in terms of $g$ only
in a range in which $dg/du$ does not vanish, namely 
$u\in [0,u^\star)$ where $u^\star$ is the smallest positive solution of $dg/du(u^\star)=0$. This 
last condition may be rewritten as 
\eqn\yu{y=F(u^\star)+ u^\star F'(u^\star)}
which has a unique solution $0<u^\star<1/27$ for any finite positive value of $y$ since
the function $F(u)+ u F'(u)$ is strictly increasing and maps $[0,1/27)$ to $[0,+\infty)$. 
This implies that $u(g)$ has a finite radius of convergence $g^\star=u^\star(y-F(u^\star))/y^2$.
Note that $d^2g/du^2 (u^\star) <0$ which implies that the singularity of $u(g)$ has of the 
form $(u^\star-u(g)) \sim (g^\star-g)^{1/2}$. Finally, we easily see that the function 
$R(u)$ in \parametric\ is regular in the range $[0,u^\star]$ (as $F(u)\leq F(u^\star)< y$ 
in this range) and $dR/du(u^\star)>0$. Hence $R(g)$ has the same radius of convergence
$g^\star$ and the same square root singularity. This corresponds to a critical exponent
$\gamma=-1/2$ (defined by $R(g^\star)-R(g) \sim (g^\star-g)^{-\gamma}$) characteristic of
the universality class of so-called {\it pure gravity}, which is also found for
quadrangulations without blocked edges.

\fig{A rooted ternary tree (left). Upon matching leaves by pairs via a set of
non-crossing arches, we obtain a tetravalent map endowed with a spanning tree, and
with a distinguished oriented edge.}{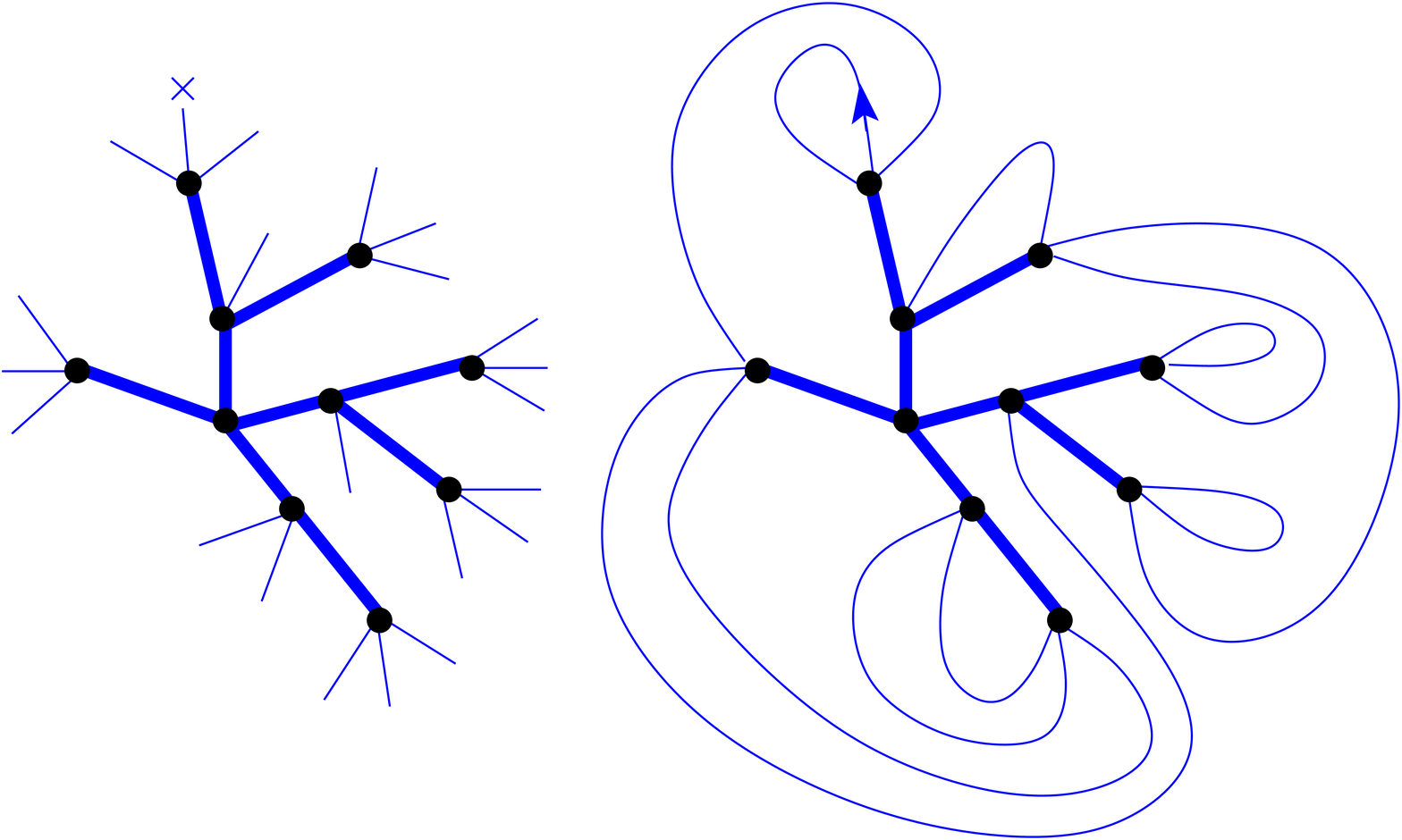}{10.cm}
\figlabel\coveringtree
Finally, a different singular behavior is obtained by taking the $y\to \infty$ limit in eq. \eqforest, keeping
$\alpha=g\, y$ fixed. In this limit, we readily see that $R=1+(1/y) F(\alpha) +{\cal O}(1/y^2)$ where
$F(\alpha)$ can now be interpreted as the generating function of maps with tetravalent vertices 
and equipped with a {\it spanning
tree}, i.e. a single connected tree passing via all the vertices of the map, and with a weight
$\alpha$ per vertex. The map has as before a distinguished face (origin face) and a distinguished
unblocked edges separating two faces at distance $m$ and $m+1$ from the origin. We simply
interpret the coefficient of $\alpha^n$ in $F(\alpha)$ as
\eqn\coefcov{{(3n)!\over n!n!(n+1)!}={(3n)!\over n!(2n+1)!}\times {(2n+2)! \over (n+1)!(n+2)!}\times
(n+2)\times {1\over 2}}
expressing, as displayed in figure \coveringtree, that our tetravalent map with a spanning tree is 
made of a rooted ternary
tree with $n$ inner vertices and $2n+2$ leaves (first factor), completed with a system of $n+1$
arches connecting these leaves (second factor) and forming the unvisited edges, with a choice of
one among $n+2$ faces (third factor). The last factor $1/2$ is because the rooting of the ternary
tree amounts to distinguish an {\it oriented} unvisited edge (that emerging from the root) while,
in $F$, we distinguish only an un-oriented edge (here all edges are of type $m\to m+1$,
irrespectively of the choice of origin).
For $\alpha \to \alpha^\star\equiv 1/27$, we now have the singular behavior $F(\alpha^\star)-F(\alpha) 
\sim (\alpha^\star-\alpha)
{\rm Log}(\alpha^\star-\alpha)$, hence a critical exponent $\gamma=-1$. The same exponent
was obtained in Ref.\DK\ for a slightly different model.

The above analysis applies straightforwardly to the case of section 3.1 by substituting
$g\to g(1+y)^2$. Among other possible direct applications of eq.\system, one may also
consider the case of arbitrary Eulerian maps with a {\it maximal number of blocked edges}.  
This case is discussed in detail in appendix B.

\newsec{Applications to maps with hard particles}

\subsec{Particles with exclusion rules on maps}

Beside the direct applications already mentioned, the main interest of Eulerian maps
with blocked edges is that they give access to the statistics of a number of quite involved
``matter" models of maps equipped with interacting particles or spins. 
More precisely, we shall focus on the quite general case of Eulerian maps with particles 
satisfying a nearest neighbor {\it $p$-exclusion rule}. This rule is expressed as follows: 
each face of the map may be occupied by a number of identical particles with the restriction 
that there be at most a total of $p$ particles on any pair of adjacent faces \HPMAT. 
For $p=1$, the exclusion rule is nothing but the 
celebrated hard-particle condition that two particles cannot occupy identical nor adjacent sites. 
In general, in addition to the usual weights $g_k$ and $\tilde g_k$ per white and black 
$k$-valent face, we attach the weight $z_i$ to each face occupied by $i$ particles (with $z_0=1$).
All these models are known to undergo various crystallization transitions to ordered phases 
in which the density of particles on black and on white faces are different (see Ref.\HPMAT\ for details). 
\fig{A piece of map (a) with blocked edges and with some maximal domain (inside the dashed line) which is not 
attainable from the origin of the map, supposed to lie outside of the figure. All potentially incoming 
edges are blocked so that no path can enter the domain. In the dual representation (b), the boundary
of the domain is a loop of length $2L$ of alternating black and white vertices, each carrying a number 
$n_i$ of particles, $i=0,\ldots,2L-1$ (and with $n_{2L}\equiv n_0$) and with every other edge blocked. 
We have $n_{2i-1}+n_{2i} \geq p+1$
for $i=1\ldots,L$, hence there exists some $i_0\in \{0,\ldots,L-1\}$ for which $n_{2i_0}
+n_{2i_0+1}\geq p+1$. The corresponding edge may be freely blocked or not, thus contributing
a weight $1+y=0$.}{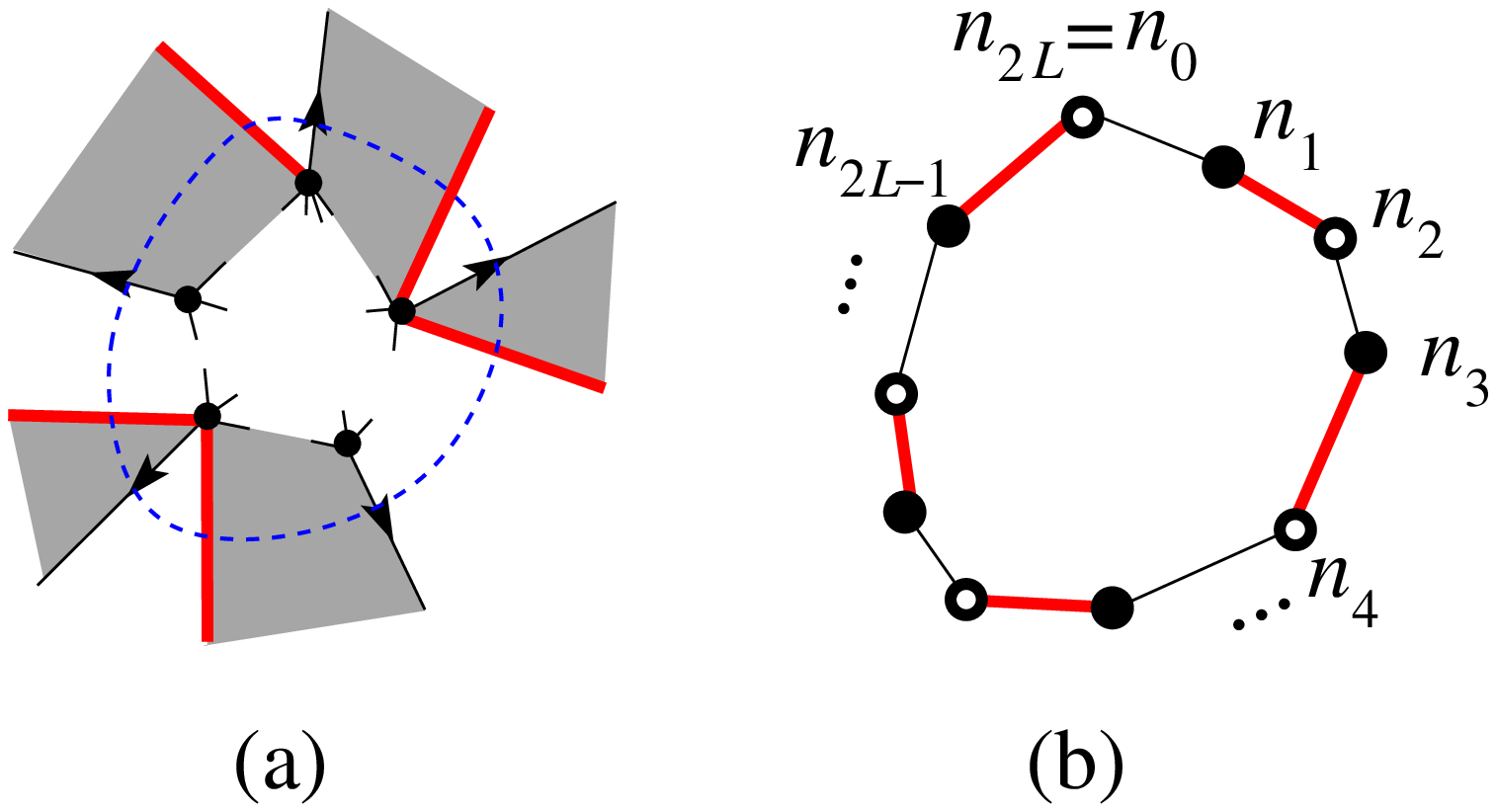}{11.cm}
\figlabel\loopy
The $p$-exclusion rule may be implemented as follows: we first start by relaxing it and demand only 
that there be at most $p$ particles on each face. We then
introduce on top of any such particle configuration a set of blocked edges on the map, 
with a weight $y$ per blocked edge, and with the two following restrictions:
\item{(i)} as before, we demand that blocked edges satisfy the global connectivity constraint; 
\item{(ii)} we furthermore restrict the configurations of blocked edges 
by allowing only for the blocking of those edges which violate the $p$-exclusion rule for the
particle configuration at hand.\par
\noindent As we already know, restriction (i) is crucial to have a mobile description of
the configurations at hand. Restriction (ii) is a new constraint but it will be completely 
transparent in the mobile formalism. 
Finally, the $p$-exclusion rule is restored by simply taking $y=-1$ in the various generating functions. 

Indeed, note first that whenever the configuration of particles satisfies the $p$-exclusion rule, 
we cannot add any blocked edge and it receives the correct weight (irrespectively of the value of $y$).
If the configuration of particles violates the $p$-exclusion rule, let us show that the contributions
of all blocked edge configurations allowed on the map add up to $0$. Clearly, if we relax 
the global connectivity constraint (i), any violated edge can be blocked or not, thus receiving
a weight $1+y=0$ so that the contribution of all blocked edge configurations vanishes trivially. 
To prove that the contribution of those configurations that satisfy (i) also vanishes, 
we therefore need only prove the vanishing of the complementary contribution of those blocked 
configurations that {\it do not} satisfy (i). To this end, assume the existence of a maximal connected
domain which is not attainable from the origin of the map and consider the boundary of the dual of this
domain (see figure \loopy-(a)). It consists of a loop of alternating black and white vertices connected
via a succession of edges dual to alternatively in- and outgoing edges of the original
Eulerian map (see figure \loopy-(b)). Clearly, all in-going edges must be blocked to 
prevent entering the domain. We simply have to show that at least one
of the outgoing edges violates the $p$-exclusion rule so that we may freely decide to block it or not,
thus creating a weight $1+y=0$ causing the whole contribution to vanish. As all in-going
edges are blocked, they violate the exclusion rule, thus the total number of particles around
the loop is larger or equal to $(p+1) L$ if the loop has length $2L$. Assuming by contradiction 
that none of the out-going edges violates the exclusion, then the total number of particles
around the loop would not exceed $p L$. This forms a contradiction.
\fig{A sample configuration of pointed Eulerian map with $4$ particles (top left) violating the
$p$-exclusion rule, here with $p=1$. We have listed all possible blocked edge configurations 
satisfying the global connectivity constraint (i) and the condition (ii) that we may block only edges
that violate the $1$-exclusion rule. Weighting each blocked edge
by $y=-1$, the contributions of the $12$ configurations add up to $0$, as it should.}{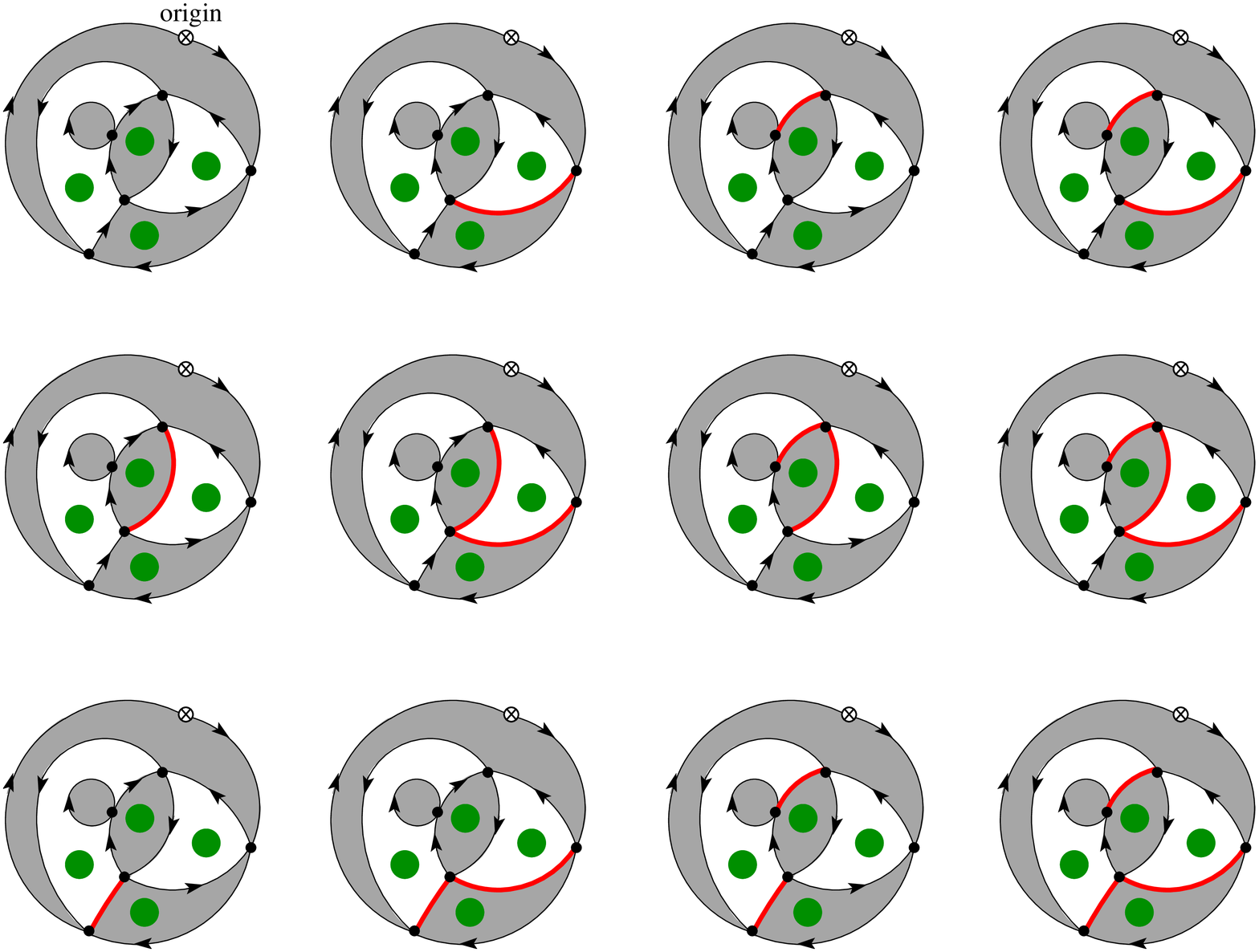}{13.cm}
\figlabel\eulermapwithhp
For illustration, we have displayed in figure \eulermapwithhp\ the allowed blocked edge
configurations satisfying (i) and (ii) above for $p=1$ and for a particular configuration
of particles. Their total contribution vanishes when $y=-1$ as the particle configuration
at hand violates the $1$-particle exclusion rule.

We conclude that the nearest-neighbor interaction of the particle model with $p$-exclusion rule
can be completely undone at the expense of introducing blocked edges on the map. 
This will allow for a straightforward derivation of recursive equations determining 
the generating function of the model.

\subsec{Generating functions}
Let us consider the generating function $G$ of Eulerian maps with particles subject to the 
$p$-exclusion rule, with a distinguished origin vertex (pointed maps) and, for convenience, with
also a {\it distinguished edge}. Thanks to the bijection of section 2, $G$ may be expressed
in terms of various generating functions for rooted sub-mobiles. More precisely, the mobiles 
that we have to consider now are unrestricted mobiles supplemented by particles living on their 
black and white unlabeled vertices, and with at most $p$ particles per vertex. Given this particle 
decoration, the only new constraint inherited form the condition (ii) of previous section is that 
a flagged edge may be marked only if the total number of particles on its two endpoints exceeds 
$p$ strictly. By a direct generalization of our notations of section 2.4, we may now introduce
the following generating functions:
\item{--} $R$ denotes as in eq.\limitgen\ the generating function for sub-mobiles rooted at a 
labeled corner, and with arbitrary (but fixed) root label;
\item{--} $B_\ell^{(i)}$ (resp. $\tilde B_\ell^{(i)}$) denotes the generating function for 
sub-mobiles rooted at some un-marked (resp. marked) flagged edge incident to a black vertex 
carrying $i$ particles and with arbitrary (but fixed) labels on the left and on the right of 
the edge differing by $\ell$ (with $\ell\geq 0$ in the un-marked case);
\item{--} finally, $W_{-\ell}^{(i)}$ (resp $\tilde W_{-\ell}^{(i)}$) denotes to the generating 
function for sub-mobiles rooted at some un-marked (resp. marked) flagged edge incident to a white vertex 
carrying $i$ particles and with arbitrary (but fixed) labels on the right and on the left
of the edge differing by $\ell$ (with again $\ell\geq 0$ in the un-marked case).\par
\noindent  By use of a straightforward generalization of eqs.\nbe\ and \be\ with $y=-1$,
the generating function $G$ is now given by:
\eqn\Gform{G=R+\sum_{i=0}^p \sum_{j=0}^p \sum_{\ell \geq 0} B_\ell^{(i)} W_{-\ell}^{(j)}
- \sum_{i=0}^{p} \sum_{j=p+1-i}^{p} \sum_\ell \tilde B_\ell^{(i)} \tilde W_{-\ell}^{(j)}.}
Collecting the generating functions of the different sub-mobiles that can be 
attached to a white (resp. black) vertex carrying $i$ particles, we are led to 
introduce $p+1$ Laurent series $Q_\circ^{(i)}(z)$ (resp. $Q_\bullet^{(i)}(z)$), 
$i=0,\ldots,p$, defined as: 
\eqn\Qcharged{\eqalign{
Q_\circ^{(i)}(z) &= {R\over z}\, +\sum_{j=0}^p \sum_{\ell\geq 0} B_{\ell}^{(j)}\, z^\ell
+\sum_{j=p+1-i}^p \sum_{\ell} \tilde B_{\ell}^{(j)}\, z^\ell \cr
Q_\bullet^{(i)}(z) &=  z\, +\sum_{j=0}^p \sum_{\ell\leq 0} W_{\ell}^{(j)}\, z^{\ell}
+\sum_{j=p+1-i}^p \sum_{\ell} \tilde W_{\ell}^{(j)}\, z^{\ell}.\cr
}}
All generating functions are now fixed by the recursive relations: 
\eqn\recurcomp{\eqalign
{ R &= {1\over 1-L}\ ,\quad L=\sum_k g_k \sum_{i=0}^p \left(Q_\circ^{(i)}(z)\right)^{k-1}|_{z^1}\cr
B_\ell^{(i)} &= z_i\ \sum_k \tilde g_k \left(Q_\bullet^{(i)}(z)\right)^{k-1}|_{z^\ell}, \quad \ell\geq 0\cr
\tilde B_\ell^{(i)} &= - z_i\ \sum_k \tilde g_k \left(Q_\bullet^{(i)}(z)\right)^{k-1}|_{z^\ell}\cr
W_\ell^{(i)} &= z_i\ \sum_k g_k \left(Q_\circ^{(i)}(z)\right)^{k-1}|_{z^\ell}, \quad \ell \leq 0\cr
\tilde W_\ell^{(i)} &= - z_i\ \sum_k g_k \left(Q_\circ^{(i)}(z)\right)^{k-1}|_{z^\ell}.\cr
}}
In particular, note that $\tilde B_\ell^{(i)}=- B_\ell^{(i)}$ and $\tilde W_{-\ell}^{(i)}=- W_{-\ell}^{(i)}$
for $\ell \geq 0$, so that eq. \Gform\ simplifies into
\eqn\Gformsimp{G=R+\sum_{i=0}^p \sum_{j=0}^{p-i} \sum_{\ell \geq 0} B_\ell^{(i)} W_{-\ell}^{(j)}
- \sum_{i=0}^{p} \sum_{j=p+1-i}^{p} \sum_{\ell <0}\tilde B_\ell^{(i)} \tilde W_{-\ell}^{(j)}}
while we may write eqs.\Qcharged\ as
\eqn\Qcharsimp{\eqalign{
Q_\circ^{(i)}(z) &= {R\over z}\, +\sum_{j=0}^{p-i} \sum_{\ell\geq 0} B_{\ell}^{(j)}\, z^\ell
+\sum_{j=p+1-i}^p \sum_{\ell< 0} \tilde B_{\ell}^{(j)}\, z^\ell \cr
Q_\bullet^{(i)}(z) &=  z\, +\sum_{j=0}^{p-i} \sum_{\ell\leq 0} W_{\ell}^{(j)}\, z^{\ell}
+\sum_{j=p+1-i}^p \sum_{\ell> 0 } \tilde W_{\ell}^{(j)}\, z^{\ell}.\cr
}}
\fig{Simplification of the crossing rule when $y=-1$. Here $i$ and $j$ denote the numbers
of particles occupying the white and black vertices. The original crossing rule (a) of section
2.2 displayed in figures \white\ and \black\ can be equivalently replaced by a new crossing rule 
(b) in which the allowed variation of labels $m,n$ is strictly correlated to the satisfaction ($m\geq n$)
or violation ($m<n$) of the $p$-exclusion rule.}{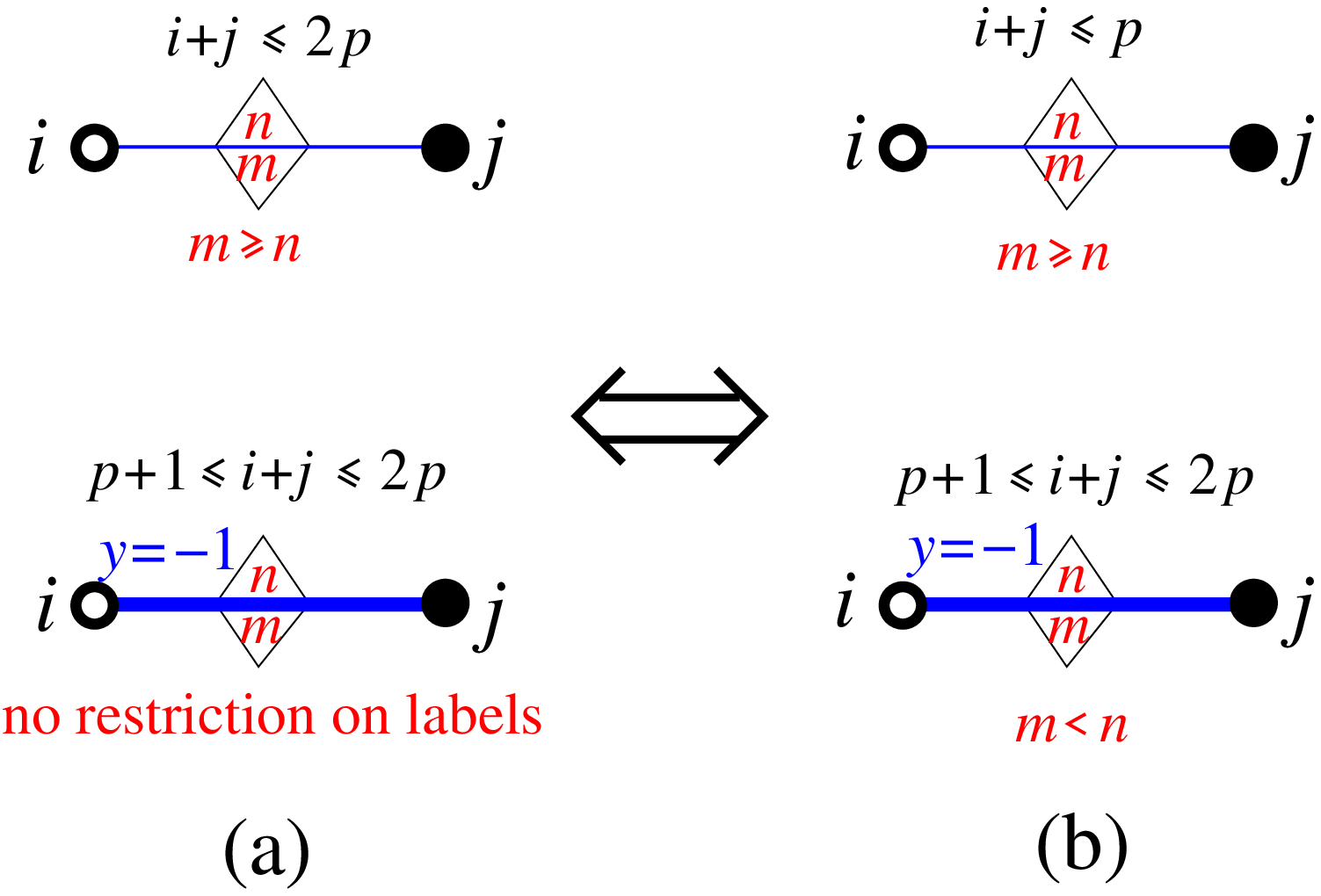}{10.cm}
\figlabel\nrule
\noindent In particular, we need only consider generating functions $\tilde B_\ell^{(j)}$
for $\ell<0$ (and $\tilde W_\ell^{(j)}$ for $\ell >0$). 
This is a direct manifestation of an obvious simplification of the crossing
rule for labels on mobiles when $y=-1$, as displayed in figure \nrule. 
Any flagged edge whose endpoints have a total number of particles larger or equal to $p+1$ 
and with increasing label clockwise around the white endpoint can be marked or not on the mobile, resulting
in a weight $1+y=0$. The mobiles that survive in practice are therefore mobiles whose flagged edges
are of the following two types only (see figure \nrule-(b)): 
\item{--} unmarked flagged edges whose endpoints have a total number of particles smaller or 
equal to $p$ and with increasing label clockwise around the white endpoint; 
\item{--} marked flagged edges whose endpoints have a total number of particles larger  
or equal to $p+1$ and with strictly decreasing label clockwise around
the white endpoint. These edges furthermore receive a weight $-1$.\par
\noindent As in section 2.4, generalizing \system, we may summarize the equations above in a closed system: 
\eqn\newsystem{\encadremath{\eqalign{
Q_\circ^{(i)}(z) 
&= {R\over z} + \sum_k \tilde g_k \sum_{j=0}^{p-i} z_j \left[ \left(Q_\bullet^{(j)}
(z)\right)^{k-1} \right]_+ - 
\sum_k \tilde g_k \sum_{j=p+1-i}^p z_j \left[ \left(Q_\bullet^{(j)}(z)\right)^{k-1}\right]_{<0}\cr
Q_\bullet^{(i)}(z) 
&= z + \sum_k g_k \sum_{j=0}^{p-i} z_j \left[ \left(Q_\circ^{(j)}(z)\right)^{k-1} 
\right]_- -
\sum_k g_k \sum_{j=p+1-i}^p z_j \left[\left(Q_\circ^{(j)}(z)\right)^{k-1}\right]_{>0}\cr
R& =1+\sum_k g_k \sum_{j=1}^p z_j \left[{\left(Q_\circ^{(j)}(z)\right)^{k-1}\over z}\right]_0  R \ ,\cr
}}}
where the notations $[\cdot]_{<0}$ and $[\cdot]_{>0}$ stand respectively
for the extraction of the negative and positive power parts in the Laurent series of $z$.

\subsec{Detailed study of Eulerian triangulations with hard particles}

As an example, let us study in detail the case of Eulerian triangulations with
hard-particles, i.e particles subject to the $1$-exclusion rule. This case
corresponds to having $g_k=\tilde g_k =g \delta_{k,3}$ (as there are as many
black and white triangles in an Eulerian triangulations, there is no need 
to introduce different weight factors for these faces).
In a dual formulation, this model is equivalent to that of so-called bicubic
maps with hard particles studied in Ref.\HPMAT\  by random matrix techniques
and in Ref.\HObipar\ by combinatorial techniques based on blossom trees.
We refer to these references for a description of the physical properties
of the model. By a simple look on eqs. \Qcharsimp\ and \recurcomp, we immediately
deduce that the only non vanishing generating functions in this case are
$B_{2}^{(0)}$, $B_{\ell}^{(1)}$ for $\ell=2,5,8$ and $\tilde B_{\ell}^{(1)}$ for $\ell=-1,-4$
together with their ``white" counterparts $W_{-2}^{(0)}$, $W_{\ell}^{(1)}$ for $\ell=-2,-5,-8$
and $\tilde W_{\ell}^{(1)}$ for $\ell=1,4$, so that we may write eq. \Qcharsimp\ as:
\eqn\Qval{\eqalign{
Q_\circ^{(0)}(z)&= {R\over z} + B_2^{(0)} z^2 +B_2^{(1)} z^2 +B_5^{(1)} z^5 +B_8^{(1)} z^8 
\cr
Q_\circ^{(1)}(z)&= {R\over z} + B_2^{(0)} z^2 +{\tilde B_{-1}^{(1)} \over z} +{\tilde B_{-4}^{(1)}
\over z^4} \cr
Q_\bullet^{(0)}(z)&= z + {W_{-2}^{(0)}\over z^2} +{W_{-2}^{(1)}\over z^2} +{W_{-5}^{(1)}\over z^5} 
+{W_{-8}^{(1)}\over z^8} \cr
Q_\bullet^{(1)}(z)&= {z} + {W_{-2}^{(0)}\over z^2} +\tilde W_{1}^{(1)} z+\tilde W_{4}^{(1)} z^4\ . \cr}}
Equations \recurcomp\ now translate into:
\eqn\recureqs{\matrix{
B_2^{(0)}=g \hfill &  W_{-2}^{(0)}=g R^2 \hfill \cr
B_2^{(1)}=g\, z_1 \left( \left(1+\tilde W_{1}^{(1)}\right)^2 +2 W_{-2}^{(0)} 
\tilde W_{4}^{(1)}\right)\hfill & 
W_{-2}^{(1)}=g\, z_1 \left( \left(R+\tilde B_{-1}^{(1)}\right)^2 +2 B_{2}^{(0)} 
\tilde B_{-4}^{(1)}\right)\hfill \cr
\tilde B_{-1}^{(1)}= - 2 g\, z_1 W_{-2}^{(0)} \left(1+\tilde W_{1}^{(1)} \right)\hfill &   
\tilde W_{1}^{(1)}= - 2 g\, z_1 B_{2}^{(0)} \left(R+\tilde B_{-1}^{(1)} \right) \hfill \cr
\tilde B_{-4}^{(1)}= - g\, z_1 \left(W_{-2}^{(0)}\right)^2 \hfill & 
\tilde W_{4}^{(1)}= - g\, z_1 \left(B_{2}^{(0)}\right)^2  \hfill \cr
}}
as well as expressions for $B_5^{(1)}$, $B_8^{(1)}$, $W_{-5}^{(1)}$ and $W_{-8}^{(1)}$ that we shall not use.
This allows to express these generating functions in terms of $R$ only, namely:
\eqn\resrec{\eqalign{
B_2^{(0)}&=g\cr
B_2^{(1)}&= { g\, z_1 (1 -2 g^4 z_1 R^2 -8 g^6 z_1^2 R^3 -8 g^8 z_1^3 R^4)
\over ( 1 +2 g^2 z_1 R)^2}\cr
\tilde B_{-1}^{(1)}&=- {2 g^2 z_1 R^2 \over 1 +2 g^2 z_1 R} \cr
\tilde B_{-4}^{(1)}&= - g^3 z_1 R^4\ ,\cr}}
while $W_{-2}^{(0)}= R^2 B_{2}^{(0)}$, $W_{-2}^{(1)}= R^2 B_{2}^{(1)}$, 
$W_{1}^{(1)}= B_{-1}^{(1)}/R$ and $W_{4}^{(1)}= B_{-4}^{(1)}/R^4$. 
Finally, the generating function $R$ itself is given by the recursion relation:
\eqn\exprR{\eqalign{R& =1+ g R \left(2 R \left(B_2^{(0)}+B_2^{(1)}\right) 
+ 2 z_1 \left(R+\tilde B_{-1}^{(1)}\right) B_{2}^{(0)}\right) \cr
&= 1+2 g^2 R^2 {1+2 z_1+2 g^2 (2 z_1+z_1^2) R +2 g^4 z_1^2 R^2 -8 g^6 z_1^3 R^3 -8 g^8 z_1^4 R^4 
\over ( 1 +2 g^2 z_1 R)^2} \ .\cr}
}
In particular, expanding $R$ in powers of $g$ leads to:
\eqn\expanR{\eqalign{R& = 1\!+\!2g^2(1\!+\!2z_1)\!+\!4g^4(2\!+\!8z_1\!+\!5z_1^2)\!+\!4g^6(10\!+
\!60z_1\!+\!89z_1^2\!+\!28z_1^3)\!+\!32g^8(7\!+\!56z_1\!\cr &+\!135z_1^2\!+\!107z_1^3\!+\!21z_1^4)\!+\!16g^{10}
(84\!+\!840z_1\!+\!2828z_1^2\!+\!3808z_1^3\!+\!1911z_1^4\!+\!264z_1^5)\cr &+\!64g^{12}(132\!+
\!1584z_1\!+\!6870z_1^2\!+\!13320z_1^3\!+\!11629z_1^4\!+\!4088z_1^5\!+\!429z_1^6) 
\!+\!{\cal O}(g^{14})\ .\cr}}
The desired generating function $G$ for maps can be expressed in terms of $R$ as:
\eqn\Gbicub{\eqalign{G&=R+B_{2}^{(0)}W_{-2}^{(0)}+ B_{2}^{(0)} W_{-2}^{(1)}+ B_{2}^{(1)} W_{-2}^{(0)}
-\tilde B_{-1}^{(1)} \tilde W_{1}^{(1)}- \tilde B_{-4}^{(1)} \tilde W_{4}^{(1)}\cr
&= R\ {
\left(1+4 g^4 z_1 R^2 +g^2 (R+6 z_1 R) -g^6 z_1^2 R^3 -20 g^8 z_1^3 R^4 -20 g^{10} z_1^4 R^5 \right)
\over ( 1 +2 g^2\, z_1 R)^2 } \ .
\cr}}
In particular, expanding $G$ in powers of $g$ leads to:
\eqn\expanG{\eqalign{G& = 1\!+\!3g^2(1\!+\!2z_1)\!+\!12 g^4(1\!+\!4z_1\!+\!2z_1^2)\!+\!15g^6(4\!+
\!24z_1\!+\!33z_1^2\!+\!8z_1^3)\!+\!48g^8(7\!+\!56z_1\!\cr &+\!130z_1^2\!+\!92z_1^3\!+\!14z_1^4)\!+\!168g^{10}
(12\!+\!120z_1\!+\!395 z_1^2\!+\!500z_1^3\!+\!220z_1^4\!+\!24z_1^5)\cr &+\!144g^{12}(88\!+
\!1056z_1\!+\!4512z_1^2\!+\!8416z_1^3\!+\!6801z_1^4\!+\!2080z_1^5\!+\!176z_1^6)
\!+\!{\cal O}(g^{14})\ ,\cr}}
which displays the generating functions of Eulerian triangulations with hard particles
and with $2,4,\ldots,12$ triangles. Note that here, the triangulations have a distinguished
origin vertex {\it and} a distinguished edge. It is more usual to consider so-called rooted
maps with a distinguished edge only. As there are $n+2$ vertices in an Eulerian triangulation
with $2n$ faces, this is done here by dividing out the coefficient of $g^{2n}$
by an overall factor $(n+2)$. 

\subsec{Ising model revisited}
\fig{An example (a) of quadrangulation with Ising up/down spins on its faces and the equivalent
Eulerian map with bi- and tetra-valent black or white faces (b). The original spin up (resp. down) 
faces translate into empty black (resp. white) tetravalent faces, while a bivalent face
is introduced between any pair of adjacent faces with equal spins to restore the bicolorability. 
The same configurations are alternatively selected from all possible Eulerian maps with
bi- and tetra-valent faces by adding a particle on each bivalent face and imposing
the $1$-particle exclusion. This indeed forces the bivalent faces to be isolated and 
to separate tetravalent faces of identical color.}{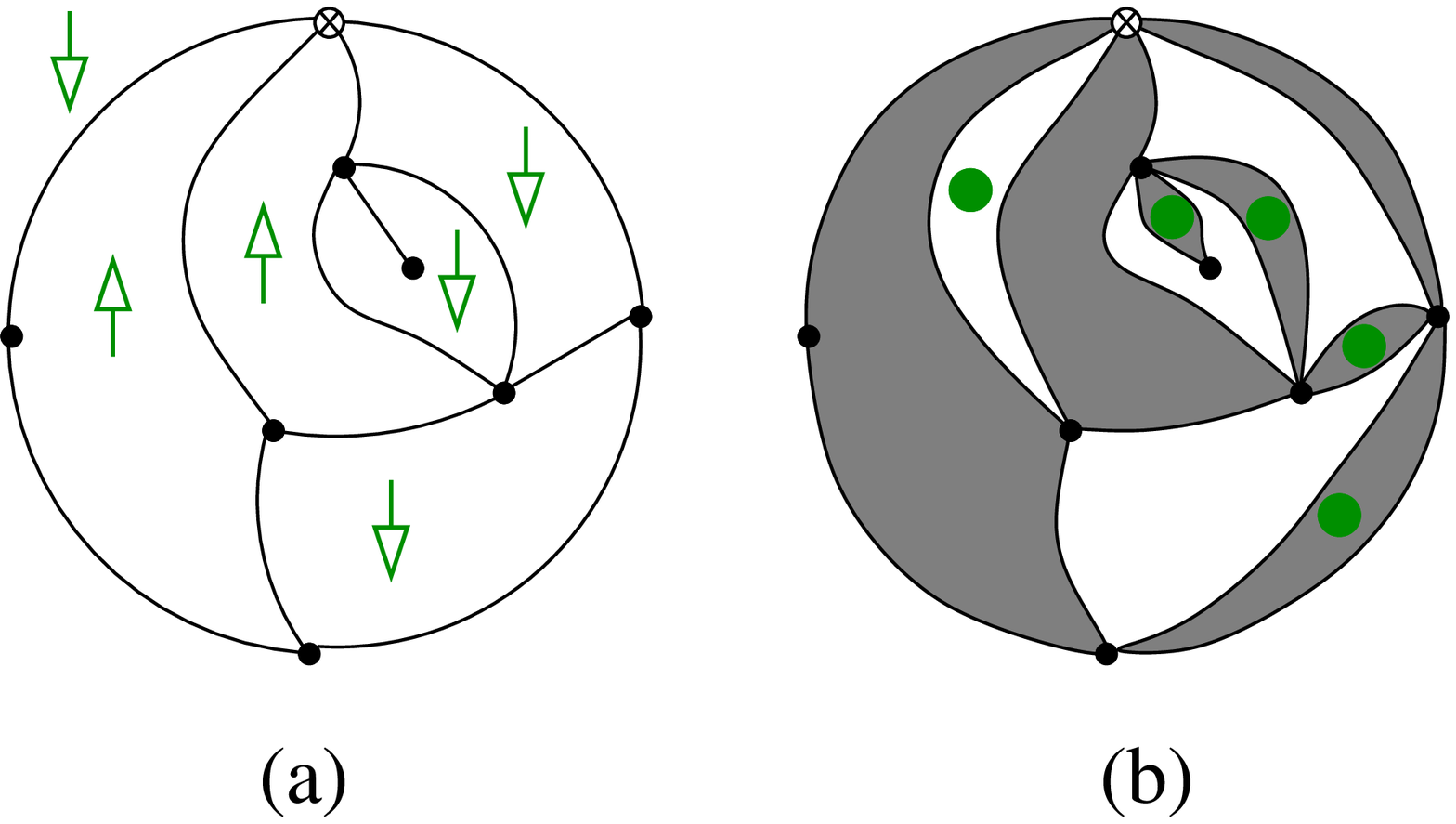}{12.cm}
\figlabel\isingviahp
As a final application, let us consider the case of quadrangulations with
Ising spins, i.e maps with tetravalent faces only and with spins up or down living 
on these faces. The Ising coupling is introduced by a weight, say $z_1$ per edge 
separating faces with the same value of the spin, while edges separating faces 
with opposite spins receive the weight $1$ instead. The model may be realized in terms of
Eulerian maps with hard-particles as follows (see figure \isingviahp): consider Eulerian maps made
of:
\item{--} black or white tetravalent faces that carry no particle;
\item{--} black or white bivalent faces that, on the contrary, are forced to carry a particle, 
weighted by $z_1$.\par  
\noindent Assume now that the particles are subject to the $1$-particle exclusion rule. 
This imposes that two bivalent faces cannot be adjacent to one another so that any bivalent 
face on the map necessarily lies between two tetravalent faces. If the bivalent face is black, 
the two neighboring tetravalent faces are necessarily white (and conversely) so that, 
by squeezing all bivalent faces into single edges, we obtain a quadrangulation with black
and white faces and with an effective weight factor $z_1$ per edge separating two faces 
of the same color. This is nothing but the desired Ising model by interpreting the black 
(resp. white) faces as carrying spins up (resp. down). The Ising model thus corresponds to
considering Eulerian maps with particles subject to the $1$-particle exclusion rule, with 
$g_k=\tilde g_k =g \delta_{k,4}+z_1 \delta_{k,2}$ and with the additional
restriction that bivalent faces must carry particles while tetravalent faces must be empty.
Again this last restriction is transparent in the mobile formalism as it simply
means that the black and white vertices of the mobile must carry a particle
or not according to their valence on the mobile. 
By a simple look on eqs. \Qcharsimp\ and \recurcomp, we immediately
deduce that now, the only non vanishing generating functions are
$B_{\ell}^{(0)}$ and $B_{\ell}^{(1)}$ for $\ell=1,3$ and $\tilde B_{\ell}^{(1)}$ for $\ell=-1,-3$,
together with their white counterparts $W_{\ell}^{(0)}$ and $W_{\ell}^{(1)}$ for $\ell=-1,-3,$
and $\tilde W_{\ell}^{(1)}$ for $\ell=1,3$. We may thus write eq. \Qcharsimp\ as:
\eqn\Qising{\eqalign{
Q_\circ^{(0)}(z)&= {R\over z} + B_1^{(0)} z +B_3^{(0)} z^3 +B_1^{(1)} z +B_3^{(1)} z^3 
\cr
Q_\circ^{(1)}(z)&= {R\over z} + B_1^{(0)} z +B_3^{(0)} z^3  +{\tilde B_{-1}^{(1)} \over z} +{\tilde B_{-3}^{(1)}
\over z^3} \cr
Q_\bullet^{(0)}(z)&= z + {W_{-1}^{(0)}\over z} +{W_{-3}^{(0)}\over z^3} +{W_{-1}^{(1)}\over z} 
+{W_{-3}^{(1)}\over z^3} \cr
Q_\bullet^{(1)}(z)&= {z} + {W_{-1}^{(0)}\over z} +{W_{-3}^{(0)}\over z^3} +\tilde W_{1}^{(1)} z+\tilde W_{3}^{(1)} z^3\ . \cr}}
Equations \recurcomp\ now translate into:
\eqn\recureqs{\matrix{
B_3^{(0)}=g\hfill &  W_{-3}^{(0)}=g R^3 \hfill \cr
B_1^{(0)}=3 g \left( W_{-1}^{(0)}+W_{-1}^{(1)}\right) \hfill &  
W_{-1}^{(0)}=3 g\, R^2 \left(B_{1}^{(0)} + B_{1}^{(1)}\right)\hfill \cr
\tilde B_{-3}^{(1)}= - z_1 W_{-3}^{(0)}\hfill & 
\tilde W_{3}^{(1)}= - z_1 B_{3}^{(0)} \hfill \cr
\tilde B_{-1}^{(1)}= -  z_1 W_{-1}^{(0)} \hfill & 
\tilde W_{1}^{(1)}= -  z_1 B_{1}^{(0)}, \hfill \cr
B_{3}^{(1)}=   z_1 \tilde W_{3}^{(1)} \hfill &
W_{-3}^{(1)}=   z_1 \tilde B_{-3}^{(1)}, \hfill \cr
B_{1}^{(1)}=   z_1 \left(1+\tilde W_{1}^{(1)}\right) \hfill &
W_{-1}^{(1)}=   z_1 \left(R+\tilde B_{-1}^{(1)}\right) \hfill \cr
}}
from which we deduce all generating functions in terms of $R$ only:
\eqn\resising{\eqalign{
B_3^{(0)}&= g\cr
B_1^{(0)}&= {3 g\, z_1 R \over 1-3 g (1-z_1^2)R}\cr
B_3^{(1)}&= - g\, z_1^2 \cr
B_1^{(1)}&= {z_1- 3 g\, z_1 R \over 1-3 g (1-z_1^2)R}\cr
\tilde B_{-3}^{(1)}&= -g\, z_1 R^3 \cr
\tilde B_{-1}^{(1)}&= -{3 g\, z_1^2 R^2\over 1-3 g (1-z_1^2)R}\cr}}
and $W_{-\ell}^{(i)}= R^\ell B_{\ell}^{(i)}$ for all the values of $\ell$ and $i$ 
mentioned above. 
As for $R$ itself, it now satisfies the recursion relation:
\eqn\exprRising{\eqalign{R& =1+ R \left(3 g R^2 \left(B_3^{(0)}+B_3^{(1)}\right) 
+ 3 g\, R \left(B_1^{(0)}+B_1^{(1)}\right)^2 +  z_1 B_{1}^{(0)}\right) \cr
&= 1+ 3g\, R^2 {2z_1^2-6 g^2 R^2 (1-z_1^2)^2 +9 g^3 R^3 (1-z_1^2)^3 +g\, R (1-4 z_1^2 +3 z_1^4)
\over ( 1 -3 g (1-z_1^2) R)^2} \ .\cr}
}
Note that this equation may be rewritten as 
\eqn\oldexp{{R\over (1-z_1^2)}=1+3 g^2 (1-z_1^2) R^3 +{z_1^2\over 1-z_1^2} {R\over  ( 1 -3 g (1-z_1^2) R)^2}\ ,}
which is precisely the form found in Refs. \BOUKA, \BMS\ and \HObipar.
The first terms of the expansion of $R$ in powers of $g$ read:
\eqn\expanRising{\eqalign{R& = 1\!+\!6gz_1^2\!+\!3g^2(1\!+\!8z_1^2\!+\!15z_1^4)\!+\!18g^3 z_1^2(11\!+
\!28z_1^2\!+\!21z_1^4)\!+\!27g^4(1\!+\!31z_1^2\!\cr &+\!229z_1^4\!+\!285z_1^6\!+\!126z_1^8)\!+\!54g^{5}z_1^2
(97\!+\!907z_1^2\!+\!2521z_1^4\!+\!1929z_1^6\!+\!594z_1^8)\cr &+\!81 g^{6}(4\!+
\!279z_1^2\!+\!4833z_1^4\!+\!19958z_1^6\!+\!30678 z_1^8\!+\!16419z_1^{10}\!+\!3861z_1^{12}) 
\!+\!{\cal O}(g^{7})\ .\cr}}
The generating function $G$ for maps is now expressed in terms of $R$ as:
\eqn\Gising{\eqalign{G&=R+B_{1}^{(0)}W_{-1}^{(0)}+ B_{3}^{(0)} W_{-3}^{(0)}+ B_{1}^{(1)} W_{-1}^{(0)}
+ B_{3}^{(1)} W_{-3}^{(0)}+ B_{1}^{(0)} W_{-1}^{(1)}+B_{3}^{(0)} W_{-3}^{(1)}\cr &
\ \ \ -\tilde B_{-1}^{(1)} \tilde W_{1}^{(1)}- \tilde B_{-3}^{(1)} \tilde W_{3}^{(1)}\cr
&= R {1\!+\!10g^2R^2(1\!-\!3z_1^2)\!-\!6gR(1\!-\!2z_1^2)\!+\!9g^4R^4(1\!-\!z_1^2)^2(1-3 z_1^2)
\!-\!6g^3R^3(1\!-\!4z_1^2\!+\!3 z_1^4)\over ( 1 -3 g (1-\, z_1)^2 R)^2 \, }
\cr}}
and expanding $G\equiv G(g,z_1)$ in powers of $g$ leads to:
\eqn\expanGising{ \eqalign{G(g, z_1)
& = 1\!+\!12 gz_1^2\!+\!4g^2(1\!+\!12z_1^2\!+\!18z_1^4)\!+\!180g^3 z_1^2(2\!+
\!5z_1^2\!+\!3z_1^4)\!+\!18g^4(2\!+\!85z_1^2\!\cr &+\!624z_1^4\!+\!693z_1^6\!+\!252z_1^8)\!+\!756g^{5}z_1^2
(12\!+\!119z_1^2\!+\!312z_1^4\!+\!207z_1^6\!+\!54z_1^8)\cr &+\!432 g^{6}(1\!+
\!91z_1^2\!+\!1642z_1^4\!+\!6681z_1^6\!+\!9450 z_1^8\!+\!4356z_1^{10}\!+\!891z_1^{12}) 
\!+\!{\cal O}(g^{7})\ .\cr}}
The coefficient of $g^n$ in $G$ is the generating function of maps with a distinguished 
origin vertex (among $n+2$) and with a distinguished edge {\it in the original Eulerian map formulation}
with bivalent faces.
This means that on the associated quadrangulations, edges separating two faces of the 
same color can be marked twice. The number of such edges is given by the power of $z_1$, so that, in the end,
coefficient of $g^n z_1^{2p}$ in \expanGising\ receives a weight $(n+2)(2n+2p)$ for all possible
markings rather than the usual factor $4n$ for rooted quadrangulations with a {\it distinguished oriented edge}
only. We can easily restore this factor and obtain the generating function $H\equiv H(g,z_1)$ 
for {\it rooted quadrangulations with Ising spins}, with expansion:
\eqn\finalising{\eqalign{H(g, z_1) &= 4 g {d\over dg} \int_0^1 d\alpha \int_0^1 d\beta
{\beta  \over \alpha} \left(G(\beta \alpha^2 g, \alpha z_1) -1\right) \cr
& = \!4 gz_1^2\!+2g^2(1\!+\!8z_1^2\!+\!9z_1^4)\!+\!108g^3 z_1^2(1\!+
\!2z_1^2\!+\!z_1^4)\!+\!12g^4(1\!+\!34z_1^2\!\cr &+\!208z_1^4\!+\!198z_1^6\!+\!63z_1^8)\!+\!216g^{5}z_1^2
(10\!+\!85z_1^2\!+\!195z_1^4\!+\!115z_1^6\!+\!27z_1^8)\cr &+\!54 g^{6}(2\!+
\!156z_1^2\!+\!2463z_1^4\!+\!8908z_1^6\!+\!11340 z_1^8\!+\!4752z_1^{10}\!+\!891z_1^{12})
\!+\!{\cal O}(g^{7})\ .\cr}}

\newsec{Conclusion}

In this paper, we have considered the combinatorics of the general class of Eulerian maps 
with blocked edges. This includes a number of matter statistical models on (Eulerian or unrestricted) maps 
whose local interactions are mediated by the choice of the blocked edge configurations.
All these models are amenable to a coding by mobiles, and we obtained the general 
recursive system of equations that determines the associated generating functions.
We gave a number of explicit direct or indirect applications. In particular,
this provides a new combinatorial approach to the problem of mutually excluding particles,
a set of models that gives access to most relevant RCFT critical points.

The general structure of the recursive equations above involve Laurent series
$Q_\bullet(z)$ and $Q_\circ(z)$ or, more generally, operators ${\bf Q}_\bullet$ 
and ${\bf Q}_\circ$ which appear as transfer operators describing the evolution
of labels around vertices of the mobile. Remarkably, these operator bare a
striking similarity with the ``position" operators standardly introduced in 
the chain-interacting multimatrix models solved by bi-orthogonal polynomial
techniques. In this framework these linear operators correspond to the multiplication
by an eigenvalue of some random matrix, usually expressed in the basis of 
orthogonal polynomials. This similarity was already observed in the simpler
mobile or blossom tree descriptions of maps without matter and seems to indicate
a deep connection between the existence of bijections with trees and the solvability
of the corresponding matrix models via orthogonal polynomials. This is yet to
be understood.

At a more practical level, the recursive equations that we face involve in general 
truncations of the Laurent series $Q_\bullet(z)$ and $Q_\circ(z)$ to {\it both} their 
positive and negative parts.  In the applications that we discussed, we have explored 
only particular cases in which either we may get rid of the truncations (section 3)
or only a finite number of coefficients in the Laurent series survive (section 4).
In the first case, this leads to simple and easily solved algebraic equations 
for $Q_\bullet(z)$ and $Q_\circ(z)$. In the second case, we end up with a closed
algebraic system for the finitely many surviving generating functions.
Note that we have not explored the most general case where infinitely many coefficients
would survive, as we do not really know how to solve the system of equations in this case.
Finding the corresponding solution could give access to more general (non necessarily algebraic) 
critical behaviors. Such behaviors have been observed for other types of statistical
models such as loop models \ON, Potts models \POTTS\ or the 6-vertex model \SIXV. The latter correspond 
to matrix models that were solved by loop equation or saddle point techniques,
as no orthogonal polynomial techniques were available. We may wonder whether our
general class of Eulerian maps with blocked edges could deal with some of these problems.

To conclude, the general framework of our approach is a modification of the natural
notion of geodesic distance on the maps, correlated with the presence of matter. This
is a quite natural principle of general relativity, here transposed at a discrete level.
In the present construction, the distance was modified by introducing blocked edges
but we may imagine other mechanisms. Note finally that our formalism allows in principle
to keep track of this modified distance (by keeping track of the labels in
generating functions), thus giving access to refined intrinsic 
geometrical properties of maps with matter. 

\bigskip
\noindent{\bf Acknowledgments:} 
The authors acknowledge support from the Geocomp project, ACI Masse de donn\'ees,
from the ENRAGE European network, MRTN-CT-2004-5616 (P.D.F. and E.G.)
from the ENIGMA European network, MRTN-CT-2004-5652 (P.D.F.) and from
the ANR program GIMP, ANR-05-BLAN-0029-01 (P.D.F.).

\appendix{A}{Combinatorial interpretation of eqs.\eqforR\ and \eqforest: 
forests on tetravalent maps as decorated even-valent maps} 

Equations \eqforR\ and \eqforest\ were obtained by use of our mobile formalism but can be given 
a simple interpretation. These equations are strongly reminiscent of that obeyed 
by the generating function 
(still denoted by $R$) of planar maps with a distinguished face and a distinguished edge and 
with vertices of arbitrary even valence $2k$, counted with weight $v_{2k}$ (see \CENSUS):
\eqn\geneeven{R=1+\sum_{k} v_{2k} {2k-1 \choose k } R^{k}\ .}
Eq.\eqforest\ is recovered by setting 
\eqn\vkforest{v_{2k}= g^{k-1} y^{k-2} {(3k-3)!\over (k-1)!(2k-1)!}}
for all $k\geq 2$ and $v_2=0$. This relation may be understood as follows: any connected component
of the forest with, say $k-1$ tetravalent inner vertices (hence $k-2$ inner edges and $2k$ leaves) 
can be contracted into a 
single $2k$-valent vertex by squeezing all its inner edges. Conversely, any $2k$-valent vertex
can be expanded into a ternary connected tree with $k-2$ edges in exactly $(3k-3)!/((k-1)!(2k-1)!)$ ways.
The enumeration of forests on maps with tetravalent vertices is therefore equivalent to that
of maps without forests but with vertices of arbitrary even valences, with the weight factor
$v_{2k}$ above simply accounting for the appropriate degeneracy factor in the squeezing process.

Concerning eq.\eqforR, it differs from eq.\eqforest\ only by an extra factor $(1+y)^{2n}$
which can be simply understood from a $1$ to $2^{2n}$ correspondence, for pointed quadrangulations 
with $2n$ edges, between allowed configurations of blocked edges where the blockings are necessarily 
in both directions and allowed
configurations of blocked edges where blockings in a single direction are also allowed. 
The correspondence is as follows: starting from a configuration with edges that are either not blocked
or blocked in both directions, we consider for each edge the orientation that leads from
its endpoint with a larger distance to the origin to that with a smaller distance (on a quadrangulation,
distances between neighbors have different parities so one must be strictly larger than the other).
For this particular orientation, we keep the blocked or unblocked nature of the original
edge but for the reverse orientation, we decide arbitrarily to block it or not. 
This leads to $2^{2n}$ configurations with the {\it same distance} on the map
as the presence or absence of these last blockings clearly does not affect this distance. 
Collecting the $y$ factors for these $2^{2n}$ configurations clearly
produces a term $(1+y)^{2n}$.  The reverse construction is as follows: starting from a configuration 
with possible blockings 
in a single direction, we first suppress all blockings oriented from a vertex at smaller distance 
from the origin to a vertex at larger distance, which leaves us with edges which are either not
blocked or blocked in one direction only (one-way edges). We then transform these one-way edges 
into edges blocked in both directions. 

\appendix{B}{Maximally blocked Eulerian maps}

\fig{A bicolored tree with as many black leaves as white leaves (left). There is a unique
matching between black and white leaves that can be realized with non-crossing arches and
is such that every face can be reached from the external face upon respecting
the canonical orientation, i.e. crossing arches with their black endpoint on the
right.}{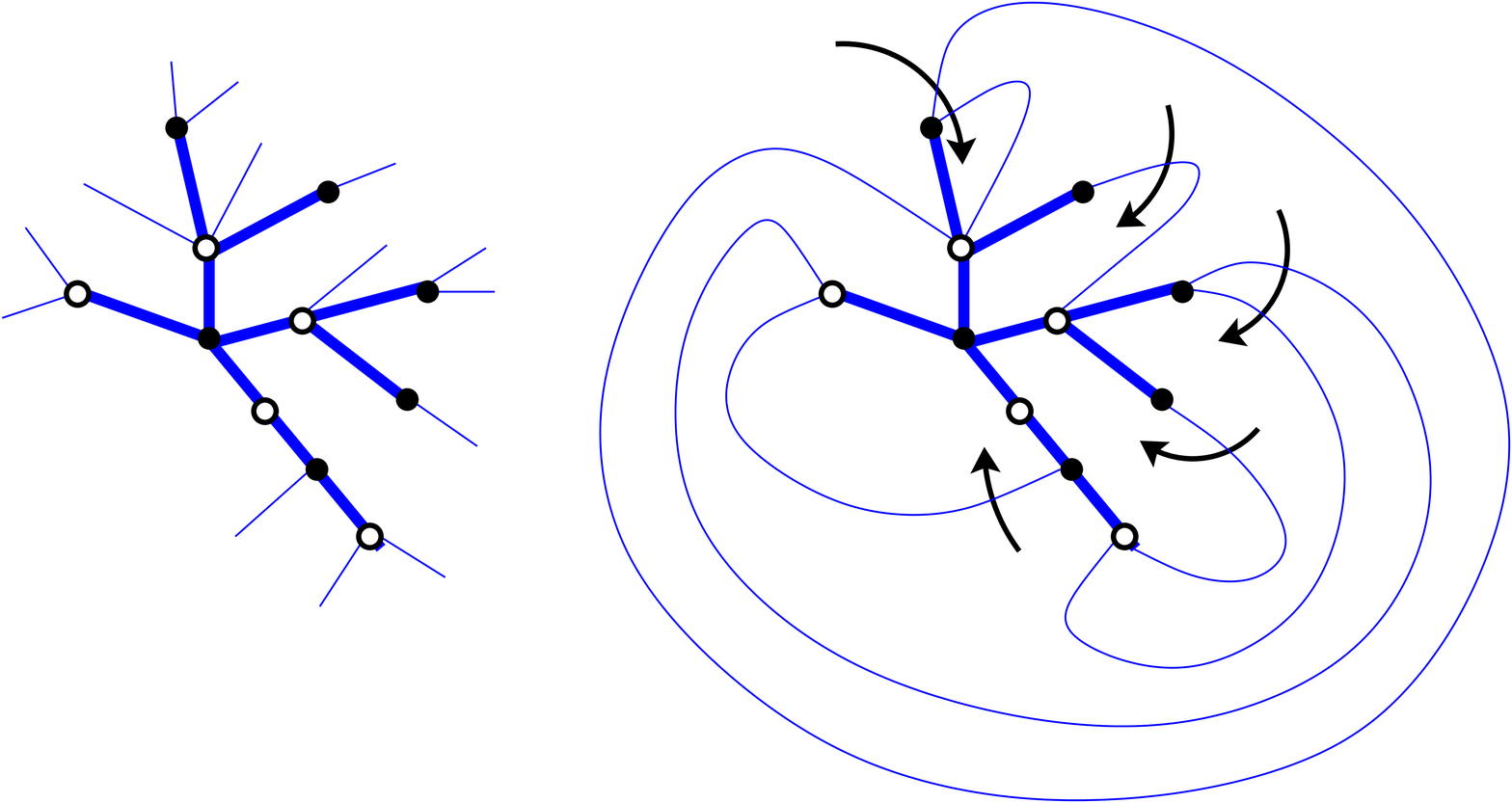}{13.cm}
\figlabel\unique
The contribution of Eulerian maps with a maximal number of blocked edges is 
obtained by taking in eqs.\system\ the limit $y\to \infty$
keeping $\alpha_k\equiv g_k\, y$ and $\tilde \alpha_k\equiv \tilde g_k\, y$ finite for all
$k$. In this case, we may write $R=1+(1/y)Z +{\cal O}(1/y^2)$ where $Z$ is the
generating function of maximally blocked Eulerian maps. From eq.\system, we have
the relations 
\eqn\maxim{\eqalign{
Z& =\sum_k \alpha_k \left[{Q_\circ^{k-1}(z)\over z}\right]_0\cr
Q_\circ(z) &= {1\over z} + 
\sum_k \tilde \alpha _k Q_\bullet^{k-1}(z)\cr
Q_\bullet(z) &= z + \sum_k \alpha_k Q_\circ^{k-1}(z)\ .\cr
}}
These equations are recognized as the recursion relations for the generating functions of 
rooted planar trees whose nodes (inner vertices) have alternating black or white color,
counted with weights $\alpha_k$ per $k$-valent white node and $\tilde \alpha_k$ per $k$-valent 
black node. More precisely, interpreting the parameter $z$ as a weight per leaf adjacent 
to a black node (hereafter referred to as black leaf), while a weight $1/z$ is attached to each leaf 
adjacent to a white node (white leaf), we 
see that $(Q_\bullet(z)/z)-1$ (resp. $z\,Q_\circ(z)-1$) is the generating function for 
such trees with a distinguished white (resp. black) leaf.
Rewriting the first line of eq.\maxim\ as
\eqn\Rmax{Z =\left[{Q_\bullet(z)\over z}-1\right]_0, }
we see that $Z$ can be interpreted as the generating function for {\it bicolored trees} with
a distinguished white leaf and constrained to have as many white leaves as black leaves. 
This should not be a surprise. Indeed, 
in a maximally blocked Eulerian map, the set of edges dual to blocked edges forms a spanning
tree of the dual map, which is naturally bicolored. Each non-blocked edge on the dual map 
connects nodes of opposite colors: cutting it into two half-edges terminating at leaves results into
a bicolored tree with as many white as black leaves. The original map has a distinguished 
non-blocked edge, which turns into a distinguish white leaf. Conversely, given a bicolored 
tree with as many white leaves as black leaves, there is a {\it unique} way to connect each black leaf
to a white one via non-intersecting arches so that each face of the resulting map can be reached
from, say the external face upon respecting the canonical orientation, i.e. crossing edges
only with the white endpoint on the left and the black endpoint on the right. This can be done
by an iterative matching algorithm: we follow the contour of the tree clockwise and match
each black leaf immediately followed by a white leaf to this latter. We then iterate the 
process by ignoring the already matched leaves until all leaves are matched (see figure \unique). 
The resulting map is clearly dual to a maximally blocked Eulerian map, with a distinguished origin vertex 
(dual to the external face) and a distinguished non-blocked edge of type $m\to m+1$ (dual to
the edge connecting the distinguished white leaf), and whose blockings satisfy the global
connectivity constraint.

\listrefs
\end